\documentclass[1 [leqno,11pt]{amsart}
\usepackage{amssymb, amsmath,amsmath,latexsym,amssymb,amsfonts,amsbsy, amsthm}

\numberwithin{equation}{section}

\allowdisplaybreaks

\numberwithin{equation}{section}

\def\sqr#1#2{{\vbox{\hrule height.#2pt
     \hbox{\vrule width.#2pt height#1pt \kern#1pt
           \vrule width.#2pt}
     \hrule height.#2pt}}}

\newtheorem{theorem}{Theorem}[section]
\newtheorem{lemma}[theorem]{Lemma}
\newtheorem{remark}[theorem]{Remark}

\newtheorem{proposition}[theorem]{Proposition}

\numberwithin{equation}{section}

\begin{document}

\title[]
{Construction of blow-up solution for 5 dimentional critical fujita type equation with different blow-up speed}

\author{Liqun Zhang}
\address{Hua Loo-Keng Key Laboratory of Mathematics, Institute of Mathematics, AMSS, and School of
Mathematical Sciences, UCAS, Beijing 100190, China 
123}
\email{lqzhang@math.ac.cn}

\author{Jianfeng Zhao}
\address{Institute of Mathematics, AMSS, UCAS, Beijing 100190, China
123}
\email{zhaojianfeng@amss.ac.cn}

\date{}
\maketitle

\begin{abstract}
We are concerned with blow-up solutions of the 5-dimensional energy critical heat equation $u_t=\Delta u + | u |^{\frac{4}{3}}u$. Our main result is to show that the existence of type \uppercase\expandafter{\romannumeral2} solutions blows up at 2 points with 2 different blow-up rates. The inner-outer gluing method has been employed.
\end{abstract}
\section{Introduction}
We consider the equation
\begin{equation}
u_t=\Delta u + | u |^{p-1}u\  \ \ \ in\in \Omega \times (0,T).
\end{equation}
Early in 1966, Fujita \cite{MR214914} started the research about equation (1.1). If we only consider the time variable, (1.1) becomes an ordinary differential equations $\dot{g}=| g |^{p-1}g$. It's solution $$g(t)=(p-1)^{-\frac{1}{p-1}}(T-t)^{-\frac{1}{p-1}}$$ blows up at time T, which is so called the type \uppercase\expandafter{\romannumeral1} blow-up. And if a blow-up solution does not satisfies this blow-up rate is called type \uppercase\expandafter{\romannumeral2} blow-up. It's known type \uppercase\expandafter{\romannumeral2} blow-up is faster than type \uppercase\expandafter{\romannumeral1} blow-up. The more result about type \uppercase\expandafter{\romannumeral1} blow-up can be seen in \cite{MR2060042,MR876989,MR1427848,MR1237055,MR3967048}.

Let $p_S=\frac{n+2}{n-2}$ denote the critical Sobolev exponent. In subcritical case, Giga and Kohn \cite{MR784476} first prove in convex domain, only type \uppercase\expandafter{\romannumeral1} blow-up occurs. In super critical case, let $p_{JL}$ denote the Joseph-Lundgreen exponent
\begin{equation*}p_{JL}=
\begin{cases}
\infty,\ \ \ \ \ \ \ \ \ \ \ \ \ \ \ \ \ \ \ \ \ 3\leq n\leq 10,\\
1+\frac{4}{n-4-2\sqrt{n-1}},\ \ \ \ \ \ \ \ \ \ n\geq11.
\end{cases}
\end{equation*}
The first example of type \uppercase\expandafter{\romannumeral2} solution is given by M.A. Herrero and Vel{\'a}zquez \cite{MR1288393,herrero1992blow}, where they constructed a positive radial solution in the case $p>p_{JL}$. Collot \cite{MR3611015} constructed a type \uppercase\expandafter{\romannumeral2} solution with the same profile later. In the case $p_S<p<p_{JL}$, Matano and Merle \cite{MR2077706} excluded type \uppercase\expandafter{\romannumeral2} blow-up in the radial case. For the case $p=p_{JL}$, Seki constructed a blow-up solution in \cite{MR3864506}. For the case $p=p_S$, Filippas, Herrero and Vel{\'a}zquez \cite{MR1843848} proved that the radial positive solution can only be a type \uppercase\expandafter{\romannumeral1} blow-up solution. By using the asymptotic matching method, they obtained formally the  possible blow-up solution with the blow-up rates given by (after a
slightly modify):
\begin{equation}\|u(\cdot,t)\|_{\infty}=
\begin{cases}
(T-t)^{-k}, \ \ \ \ \ \ \ \ \ \ \ \ \ \ \ \ \ \ \ \ \ \ \ \ \ n=3,\\
(T-t)^{-k}|\ln(T-t)|^{\frac{2k}{2k-1}},\ \ \ \ \ \ n=4,\\
(T-t)^{-3k}, \ \ \ \ \ \ \ \ \ \ \ \ \ \ \ \ \ \ \ \ \ \ \ \ n=5,\\
(T-t)^{\frac{5}{2}}|\ln(T-t)|^{\frac{15}{4}},\ \ \ \ \ \  \ \ \ \ n=6,
\end{cases}
\end{equation}
where $k=1,2,\cdots$. When $n\geq7$, Collot, Merle and Rapha{\"e}l \cite{MR3623259} proved the type \uppercase\expandafter{\romannumeral2} blow-up solution can't be around the ground state. For n=4 type \uppercase\expandafter{\romannumeral2} solutions has been constructed by Schweyer \cite{MR2990063} in 2012, In the case $3\leq n\leq 6$, type \uppercase\expandafter{\romannumeral2} solutions has been constructed by Del Pino, Musso, Wei and Zhou \cite{del2020type}, they have construct a multi point blow-up in 2019. The case n=5 is first construct by Del Pino, Musso, Wei \cite{MR3952701} in 2018, with the blow up rate of k=1 in (1.2). Junichi Harada \cite{MR4072807} finished the construction of higher blow-up rate for n=5 in 2019. And Junichi Harada \cite{harada2020type} also proved the case of n=6. Very recently Del Pino, Musso, Wei, Zhang and Zhou \cite{deltype} finished the proof of n=3.

In the above mentioned construction of the critical case, except the construction of Schweyer, they all employed the so called inner-outer gluing method. This technique is not only applied to the construction of finite time blow-up solution for the critical case, but also to super critical case and infinity time blow-up \cite{MR4046015,MR4047646}. Beside Fujita equation, D{\'a}vila, Del Pino and Wei also construct a blow-up solution for harmonic map equation \cite{MR4054257}.

\subsection{Main result}
For 5 dimensional energy critical equation, our main result is
\begin{theorem}
Let $\Omega=\mathbb{R}^n$, n = 5, $p=\frac{7}{3}$ be the critical exponent. For any given two points 0, q in $\mathbb{R}^5$, then for any sufficiently small $T>0$, there exists an initial data $u_0$, such that the solution of equation (1.1) blows up at time T, where the blow-up points are at the given 0 and q with different blow-up rate. Moreover, the main order of the solution is as
\begin{equation*}
u(x,t)=\sum_{j=1,2} U_{\lambda_j (t),\xi_j (t)}(x)-Z_1 \eta_1-Z_2+\theta(x,t),
\end{equation*}
when $t\rightarrow T$
\begin{equation*}
\lambda_j (t)\rightarrow 0,\ \ \ \ \xi_1 (t)\rightarrow 0,\ \ \ \ \xi_2 (t)\rightarrow q,
\end{equation*}
and
$\lambda_1(t)\sim C_1(T-t)^4$,
$\lambda_2(t)\sim C_2(T-t)^2$, with
\begin{equation*}
\theta(x,t) \in L^\infty(\mathbb{R}^5\times(0,T)),
\end{equation*}
\end{theorem}
The definition of $Z_1$, $Z_2$ and $\eta_1$ are left to section 2.
\begin{remark}
In our construction, actually we can construction any k points blow-up with one point has blow-up rate $(T-t)^{-6}$ and the others have blow-up rate $(T-t)^{-3}$.
\end{remark}

\section{Basic facts in the construction}\label{sec:intro}
As we all known, in critical case, the steady solution of equation (1.1) is unique, up to scalings and translations.
\begin{equation}
U(x)=\alpha_n({\frac{1}{1+|x|^2}})^{\frac{n-2}{2}},\alpha_n=(n(n-2))^{\frac{n-2}{4}},
\end{equation}
which is called Talenti-Aubin steady state \cite{MR463908}.
\subsection{The Construction of the approximated solutions}
\ \\
Fixed a point $q\in \mathbb{R}^5$, we hope to find a solution of this form:
\begin{equation}
u(x,t)=\sum_{j=1,2} U_{\lambda_j (t),\xi_j (t)}-Z_1 \eta_1-Z_2+\theta(x,t),
\end{equation}
where
\begin{equation*}
U_{\lambda_j (t),\xi_j (t)}={\lambda_j (t)}^{-\frac{n-2}{2}}U(\frac{x-\xi_j (t)}{\lambda_j (t)}).
\end{equation*}
$\lambda_j (t)$ are the scaling parameters, $\xi_j (t)$ are the translation parameters, translation is small compare with scaling, that is, $\xi_j (t)=o(\lambda_j (t))$ .\\
$Z_1:=M(T-t) e_1$, where $M$ is a small constant, $e_1=1-\frac{|z|^2}{10}$, $z=\frac{x}{\sqrt{T-t}}$ is the self similar variable. Here $Z_1$ is actually a solution of heat equation. Since the solution does not decay at the infinity,  we need to cut it off by multiplying the cut off function $\eta_1=\eta(\frac{x}{(T-t)^{\frac{1}{8}}})$, where $\eta(x)=1$ if $x\leq1$, and $\eta(x)=0$ if $x\geq2$ and it's smooth in the interval $(1,2)$.\\
We let $Z_2$ be the unique solution of
\begin{equation}\begin{cases}
\partial_tZ_2=\Delta Z_2 \ \ \ \ (x,t) \in \mathbb{R}^5 \times (0,\infty),\\
Z_2(\cdot,0)=Z_{2,0} \ \ \ \ x \in \mathbb{R}^5.
\end{cases}\end{equation}
where $Z_{2,0}$ is an odd function with compact support. $Z_{2,0}$ has a small $L^\infty$ norm. It satisfies $Z_{2,0}(q)>\frac{\|Z_2\|_\infty}{2}$ in the interval [0,T]. We can obtain it by use the continue property of heat equation. By uniqueness, it is easy to see $Z_2$ is an odd function.

\subsection{Inner-outer gluing}\ \\
We write the remainder in the form
\begin{equation}
\theta(x,t)=\lambda_1^{-\frac{n-2}{2}}\phi_1(y_1,t)\eta_R(y_1)+\lambda_2^{-\frac{n-2}{2}}\phi_2(y_1,t)\eta_R(y_1)+\psi(x,t),
\end{equation}
where $y_1=\frac{x-\xi_1}{\lambda_1}$. We denote $\eta_R(y_1)=\eta_{R,1}$, and also $\eta_R(y_2)=\eta_{R,2}$ later.\\
We require $\psi$ be a continuous function and has the following estimate
\begin{equation}
\begin{split}
|\psi|\leq\begin{cases}\delta_0 (T-t)(1+|z|^4)\ &\ |z|\leq(T-t)^{-\frac{1}{4}},\\
\frac{\delta_0}{1+|x|^2}\ &\ |z|>(T-t)^{-\frac{1}{4}},
\end{cases}
\end{split}
\end{equation}
where
\begin{equation*}\delta_0=\frac{\|Z_2\|_\infty}{10}.\end{equation*}
Let us set the error function\\
\begin{equation*}
S(u)=-u_t+\Delta u+|u|^{p-1}u.
\end{equation*}
Denote
\begin{equation*}
y_i=\frac{x-\xi_i}{\lambda_i},
\end{equation*}
\begin{equation}
E_i(y_i,t)=\lambda_i\dot{\lambda}_i[(\frac{3}{2}U(y_i)+y_i\cdot \nabla U(y_i))]+\lambda_i\dot{\xi}_i\cdot\nabla U(y_i).
\end{equation}
We compute $S(u)$
\begin{equation*}
\begin{split}
&S(\sum_{i=1,2} U_{\lambda_i (t),\xi_i (t)}-Z_1 \eta_1-Z_2+\theta(x,t))\\
&=-\theta_t+\Delta \theta+\sum_{i=1,2} pU_{\lambda_i ,\xi_i }^{p-1}(\theta(x,t)-Z_1 \eta_1-Z_2)+\sum_{i=1,2} \lambda_i^{-\frac{7}{2}}E_i+N\\
&=\sum_{i=1,2} \eta_{R,i} \lambda_i^{-\frac{7}{2}}(-\lambda_i^2\partial_t\phi_{i}+\Delta_y\phi_i+pU^{p-1}[\phi_i+\lambda_i^{\frac{3}{2}}(-Z_1 \eta_1-Z_2+\psi)]+E_i)\\
&-\psi_t+\Delta_x\psi+\sum_{i=1,2}p\lambda_i^{-2}(1-\eta_{R,i})U^{p-1}(-Z_1 \eta_1-Z_2+\psi)
+\sum_{i=1,2} A[\phi_i]+\sum_{i=1,2} B[\phi_i]\\
&+\sum_{i=1,2}\lambda_i^{-\frac{7}{2}}E_i(1-\eta_{R,i})+N.
\end{split}
\end{equation*}
Where
$$A[\phi_i]:=\lambda_i^{-\frac{7}{2}}[\Delta_{y_i}\eta_{R,i}\phi_i+2\nabla_{y_i}\eta_{R,i}\nabla_{y_i}\phi_i],$$

\begin{equation*}
\begin{split}
B[\phi_i]:=&\lambda_i^{-\frac{5}{2}}[\dot{\lambda}_i({y_i}\cdot\nabla_{y_i}\phi_i+\frac{3}{2}\phi_i)\eta_{R,i}+\dot{\xi}_i\cdot\nabla_{y_i}\phi_i\eta_{R,i}+(\dot{\lambda}_i{y_i}\cdot\nabla_{y_i}\eta_{R,i}+\dot{\xi}_i\cdot\nabla_{y_i}\eta_{R,i})\phi_i],\\
N:=&|U_{\lambda_1,\xi_1}+U_{\lambda_2,\xi_2}+\theta-Z_1 \eta_1-Z_2|^{p-1}(U_{\lambda_1,\xi_1}+U_{\lambda_2,\xi_2}+\theta-Z_1 \eta_1\\
&-Z_2)-U_{\lambda_1,\xi_1}^p-U_{\lambda_2,\xi_2}^p-pU_{\lambda_1,\xi_1}^{p-1}(\theta-Z_1 \eta_1-Z_2)-pU_{\lambda_2,\xi_2}^{p-1}(\theta-Z_1 \eta_1-Z_2)\\
&+\frac{\partial \eta_1}{\partial t}Z_1+\nabla_x Z_1\nabla \eta_1-\Delta \eta_1 Z_1.
\end{split}
\end{equation*}
We shall find a solution of (1.1) if we find a pair $(\phi_1(y,t),\phi_2(y,t),\psi(x,t))$ solves the following system of equations
\begin{equation}\begin{cases}
\lambda_1^2\partial_t\phi_{1}=\Delta_{y_1}\phi_1+pU^{p-1}(y_1)+H_1(\psi,\lambda_1,\lambda_2,\xi_1,\xi_2)\ \ \ (y_1,t)\in\ B_{2R}(0)\times(0,T),\\
\phi_{1}(y,0)=\phi_{1,0} \ \ \ \ \ \ \ \ \ \ \ \ \ \ \ \ \ \ \ \ \ \ \ \ \ \ \ \ \ \ \ \ \ \ \ \ \ \ \ \ \ \ \ \ \ \ \ \ \ \ \ \ \ \ \ y_1\in B_{2R}.
\end{cases}
\end{equation}

\begin{equation}\begin{cases}
\lambda_2^2\partial_t\phi_{2}=\Delta_{y_2}\phi_2+pU^{p-1}(y_2)+H_2(\psi,\lambda_1,\lambda_2,\xi_1,\xi_2)\ \ \ (y_2,t)\in\ B_{2R}(0)\times(0,T),\\
\phi_{2}(y,0)=\phi_{2,0}\ \ \ \ \ \ \ \ \ \ \ \ \ \ \ \ \ \ \ \ \ \ \ \ \ \ \ \ \ \ \ \ \ \ \ \ \ \ \ \ \ \ \ \ \ \ \ \ \ \ \ \ \ \ \ y_2\in B_{2R}.
\end{cases}
\end{equation}

\begin{equation}\begin{cases}
\psi_t=\Delta_x\psi+G(\phi,\psi,\lambda_1,\lambda_2,\xi_1,\xi_2)\quad (x,t)\in \mathbb{R}^5\times (0,T),\\
\psi(\cdot,0)=\psi_0\ \ \ \ \ \ \ \ \ \ \ \ \ \ \ \ \ \ \ \ \ \ \ \ \ \ \ \ \ \ \ \ \ \ \ \ \ \ \ \  x\in \mathbb{R}^5.
\end{cases}
\end{equation}
Where for $i=1,2$
\begin{equation}
H_i(\psi,\lambda_i,\xi_i):=\lambda_i^{\frac{3}{2}}pU(y_i)^{p-1}(-Z_1(\xi_i+\lambda_i y_i)-Z_2(\xi_i+\lambda_i y_i)+\psi(\xi_i+\lambda_i y_i))+E_i(y_i,t),
\end{equation}
and
\begin{equation}
\begin{split}
G(\phi_1,\phi_2,\psi,\lambda_1,\lambda_2,\xi_1,\xi_2):=&\sum_{i=1,2} p\lambda_i^{-2}(1-\eta_{R,i})U(y_i)^{p-1}(-Z_1\eta_1-Z_2+\psi)\\
&+\sum_{i=1,2} A[\phi_i]+\sum_{i=1,2} B[\phi_i]+\sum_{i=1,2} \lambda_i^{-\frac{7}{2}}E_i(1-\eta_{R,i})+N.
\end{split}\end{equation}
\subsection{A result about the linearized equation}
To deal with the inner problemma, we consider the corresponding linear problem of (1.1).
\begin{equation}
\begin{split}
\begin{cases}
\lambda^2\phi_t&=\Delta_y\phi+pU(y)^{p-1}\phi+h(y,t)\ \ \ \ (y,t)\in\ B_{2R}\times[0,T),\\
\phi (y,0)&=lW_0(y)\ \ \ \ \ \ \ \ \ \ \ \ \ \ \ \ \ \ \ \ \ \ \ \ \ \ \ \ \ \ \ \ \ \ y\in B_{2R},
\end{cases}
\end{split}
\end{equation}
where $l$ is a constant be choosing later, $W_0$ is the eigenfunction of $L_0=\Delta +pU^{p-1}$ with the negative eigenvalue.\\
It is known that
\begin{equation*}
W_0(y)\sim |y|^{-2}e^{-\sqrt{\lambda_0}|y|},\ \ \ \ as\ y\rightarrow \infty,
\end{equation*}
where $\lambda_0$ is the only negative eigenvalue of $L_0$.\\
And $h(y,t)$ satisfies the orthogonal condition
\begin{equation}
\int_{B_{2R}}h(y,t)W_i(y)dy=0\ \ \ \ i=1,\cdots,n+1,\ \ \ \ t\in\ [0,T).
\end{equation}
We define
\begin{equation*}
\|h\|_{a,\nu}:=\sup_{y\in B_2R,t\in[0,T)}\lambda_0^{\nu}(1+|y|^a)|h(y,t)|,
\end{equation*}
and
\begin{equation*}
\|\phi\|_{*a,\nu}:=\sup_{y\in B_2R,t\in[0,T)}\lambda_0^{\nu}\frac{(1+|y|^6)}{R^{6-a}}|\phi (y,t)|.
\end{equation*}
Actually we expect solution $\phi$ in the space where the $\|\cdot\|_{a,\nu}$ is finite. However, we can only know the estimate due to Del Pino, Musso and Wei, there they proved $\phi$ is in $\|\cdot\|_{*a,\nu}$. We copy the lemma of Del Pino, Musso and Wei here
\begin{lemma}
There exists constant $C_3>0$, for all sufficient large $R >0$, if $\|h\|_{2+a,\nu}<+\infty$, and h satisfies orthogonal condition, then there exists linear operator L satisfies
\begin{equation*}
\phi=L^{in}[h],\ \ \ \ l=l[h],
\end{equation*}
fulfill the equation, and l[h] satisfies
\begin{equation*}
|l[h]|+\|(1+y)\nabla_y\phi\|_{*a,\nu}+\|\phi\|_{* a,\nu}\leq C_3\|h\|_{2+a,\nu}.
\end{equation*}
\end{lemma}
\subsection{Local properties of the heat equation}\ \\
\noindent In order to solve the outer problem, we shall first list some known results for the completeness of the paper. We will not restrict in dimension $n=5$. Consider the heat equation
\begin{equation*}
\begin{split}
u_t=\Delta u\ \ \ \ (x,t)\in \mathbb{R}^n\times(0,\infty).
\end{split}\end{equation*}
We make use of self-similar variables
\begin{equation}
w(z,\tau)=u(x,t),\ \ \ z=\frac{x}{\sqrt{T-t}},\ \ \ T-t=e^{-\tau}.
\end{equation}
The function $w(z,\tau)$ solves
\begin{equation*}
\begin{split}
w_\tau=A_z w\ \ \ \ (z,\tau)\in \mathbb{R}^n\times(0,\infty),
\end{split}\end{equation*}
where $A_z=\Delta_z-\frac{z}{2}\cdot\nabla$. We define the weighted $L^2$ space\\
\begin{equation*}
L_\rho^2(\mathbb{R}^n):=\{f\in L_{loc}^2(R^n),\|f\|_\rho<\infty\},\ \|f\|_\rho^2=\int_{R^n}f^2(z)\rho(z)dz,\
\end{equation*}
where $\rho(z)=e^{-\frac{|z|^2}{4}}$.
And inner product is denoted by
\begin{equation*}
(f_1,f_2)_\rho:=\int_{R^n}f_1(z)f_2(z)\rho(z)dz.
\end{equation*}
Consider the eigenvalue problem of $A_z$ in $L_\rho^2(\mathbb{R}^n)$. Let $e_\alpha$ be the eigenfunction
\begin{equation*}
-A_ze_\alpha=\lambda_\alpha e_\alpha,
\end{equation*}
where $\alpha$ is multi-index, $\alpha=(\alpha_1,\cdots,\alpha_n)\in\ \mathbb{Z}^n$, corresponding eigenvalue is
\begin{equation}
\lambda_\alpha=\frac{|\alpha|}{2},\ \ \ \ |\alpha|=\alpha_1+\cdots+\alpha_n.
\end{equation}
The eigenfunction
\begin{equation*}
e_\alpha=\prod e_{\alpha_i},\ \ \ \ e_{\alpha_i}=H_{\alpha_i},
\end{equation*}
where $H_{\alpha_i}$ is the $\alpha_i$-th Hermite polynomial. Then $(T-t)^{\lambda_\alpha}e_\alpha$ is a solution of heat equation.
We denote $e^{A_z}f_0$ to be the solution with initial data $f_0$, by Duhamel's principle
\begin{equation}
e^{A_z}f_0=\frac{c_n}{(T-e^{-\tau})^{\frac{n}{2}}}\int_{R^n}e^{-\frac{|e^{-\frac{\tau}{2}}z-\zeta|^2}{4(T-e^{-\tau})}}f_0(\zeta)d\zeta.
\end{equation}
By using this formula, It is easy to obtain the following estimate.

\begin{lemma}
There exists a constant $C=C_n>0$, such that
\begin{equation}
|e^{A_z (\tau-\tau_0)}f_0|\leq C_n \frac{e^{\frac{e^{-(\tau-\tau_0)|z|^2}}{4(2-T+e^{-\tau})}}}{(T-e^{-\tau})^{\frac{n}{4}}}\|f_0\|_\rho.\ \ \ \ \ \ \ (z,\tau)\in \mathbb{R}^n\times(0,\infty)
\end{equation}
\end{lemma}
\noindent The proof of this lemma can be seen in \cite{MR3673663}.

\begin{lemma}
For any $l\in \mathbb{Z}^+$, it has
\begin{equation}
|e^{A_z (\tau-\tau_0)}|z|^{2l}|\leq C_l (1+e^{-l(\tau-\tau_0)}|z|^{2l}).\ \ \ \ \ \ \ (z,\tau)\in \mathbb{R}^n\times(0,\infty)
\end{equation}
\end{lemma}
\noindent The proof of this lemma can be seen in \cite{MR3484968}.\\
The last one is an estimate of heat equation in different scaling variable.

\begin{lemma}Consider nonhomogeneous heat equation
\begin{equation}
\begin{cases}
\partial_t\Phi-\Delta\Phi=g\ \ \ \ \ (x,t)\in \mathbb{R}^n\times(0,T),\\
\Phi\mid_{t=0}=0 \ \ \ \ \ \ \ \ \ \ \ x\in \mathbb{R}^n.
\end{cases}
\end{equation}
Suppose \begin{equation*}|g|\leq\frac{1}{\lambda^2(1+|y|^{2+a})},\end{equation*} then there exist constants $C_1$, ${\gamma}>0$ such that
\begin{equation*}|\Phi|\leq C_1(T^{\gamma}+\frac{1}{1+|y|^a}).\end{equation*}
\end{lemma}
\noindent The proof of this lemma can be seen the lemma 4.2 in \cite{MR3952701}.

\section{The parameters $\lambda_i$ and $\xi_i$}\label{chap:3}
In this section we shall solve the scaling parameter $\lambda_i$ and translation parameter $\xi_i$. They essentially determine the blow up rates.
\subsection{Formal derivation of $\lambda_i$ and $\xi_i$}\ \\
We formally drive approximate equation of $\lambda_i$ and $\xi_i$, and determine the first order of the parameter.

\noindent Recall the inner problem
\begin{equation*}
\begin{split}
\lambda_i^2\partial_t\phi_{i}=\Delta_{y_i}\phi_i+pU^{p-1}(y_i)+H_i(\psi,\lambda_i,\xi_i)\ \ \ (y_i,t)\in B_{2R}\times(0,T),
\end{split}
\end{equation*}
and the notation
\begin{equation*}
y_i=\frac{x-\xi_i}{\lambda_i},
\end{equation*}
where $i=1,2$ and
\begin{equation*}
H_i(\psi,\lambda,\xi):=\lambda_i^{\frac{3}{2}}pU(y_i)^{p-1}(-Z_1(\xi_i+\lambda_i y_i)-Z_2(\xi_i+\lambda_i y_i)+\psi(\xi_i+\lambda y_i))+E_i(y_i,t).
\end{equation*}
Because it is a type \uppercase\expandafter{\romannumeral2} blow-up, then\\
\begin{equation*}
|\lambda_i|<<\sqrt{T-t}.
\end{equation*}
Therefore it is reasonable for us to ignore $\lambda_i^2\partial_t\phi_{i}$, the equation becomes an elliptic equation
\begin{equation}
\Delta_{y_i}\phi_i+pU^{p-1}(y_i)+H_i(\psi,\lambda_i,\xi_i)=0.
\end{equation}
From Fredholm theorem, the solvability condition of this equation is the nonhomogeneous term orthogonal to the kernel of the operator.\\We let
\begin{equation*}L_0u=\Delta u+pU^{p-1}u.\end{equation*}
The kernels of $L_0$ are
\begin{equation}
W_j(x)=\frac{\partial U(x)}{\partial x}\ \ \ \ j=1,2,3,4,5,
\end{equation}
and
\begin{equation}
W_6(x)=\frac{\partial U_{\lambda}(x)}{\partial \lambda}\mid_{\lambda=1}=(\frac{3}{2}U(x)+x\cdot \nabla U(x)).
\end{equation}
The orthogonal conditions are approximate as
\begin{equation}
\int_{R^5}H_i(y_i,t)W_j(y_i)dy=0 \ \ \ \ j=1,\ldots,6,\quad i=1,2.
\end{equation}

Because $Z_2$ is an odd function, then $Z_2(0,t)=0$. In $y_1\in B_{2R}$, by the mean value theorem, $|Z_2(\lambda_1 y_1,t)| \leq 2R|\lambda_1|$ is a higher order term which could be omitted.  We require remainder $\psi$ can not make a big influence and also be omitted.
Notice that
$$Z_1=M(T-t)(1-\frac{|z|^2}{10}),$$ in the region $B_{2R}$ it is approximate as $M(T-t)$.\\
And from (2.5)
\begin{equation*}
E_i(y_i,t)=\lambda_i\dot{\lambda}_i[(\frac{3}{2}U(y_i)+y_i\cdot \nabla U(y_i))]+\lambda_i\dot{\xi}_i\cdot\nabla U(y_i).
\end{equation*}
Also
$$W_6(y_i)=(\frac{3}{2}U(y_i)+y_i\cdot \nabla U(y_i))$$ is an even function, while $W_j(y_i)=\frac{\partial U(y_i)}{\partial y_i}$ is odd in variable $y_i$.\\
So  $$\int_{R^5}H_1(y_1,t)W_6(y_1)dy=0$$ can approximately be written as
\begin{equation*}
\dot{\lambda}_1\lambda_1\int_{R^5}W_{6}^2-M(T-t)\lambda_1^{\frac{3}{2}}\int_{R^5}pU(y_1)^{p-1}W_{6}(y_1)dy_1=0.
\end{equation*}
Let
\begin{equation*}
M\frac{\int_{R^5}pU(y_1)^{p-1}W_{6}(y_1)dy_1}{\int_{R^5}W_{6}^2}=-\frac{3M}{2}\frac{\int_{R^5}U(y_1)^{p}dy_1}{\int_{R^5}W_{6}^2}:=-\kappa_1.
\end{equation*}
Then the main order of $\lambda_1$ is
\begin{equation*}
\lambda_{1,0}=\frac{\kappa_1^2 (T-t)^4}{16}.
\end{equation*}

Next we deal with $\lambda_2$, because $Z_1\eta_1$ is supported around 0, but the integration region is around q, so we can omit $Z_1\eta_1$. And we also omit the remainder $\psi$. $Z_2$ can be approximated as $Z_{2,0}(q)$ around $q$, then $$\int_{R^5}H_2(y_2,t)W_6(y_2)dy=0$$ can be approximated by
\begin{equation*}
\dot{\lambda}_2\lambda_2\int_{R^5}W_{6}^2-\lambda_2^{\frac{3}{2}}\int_{R^5}pU(y_2)^{p-1}W_{6}(y_2)Z_{2,0}(q)dy_2=0.
\end{equation*}
Let
\begin{equation*}
\frac{\int_{R^5}pU(y_2)^{p-1}W_{6}(y_2)Z_{2,0}(q)dy_2}{\int_{R^5}W_{6}^2}=-Z_{2,0}(q)\frac{\int_{R^5}U(y_2)^pdy_2}{\int_{R^5}W_{6}^2}=-\kappa_2.
\end{equation*}
Then the main order of $\lambda_2$ is
\begin{equation*}
\lambda_2=\frac{\kappa_2^2 (T-t)^2}{4}.
\end{equation*}
Similarly about $$\int_{R^5}H_i(y_i,t)W_j(y_2)dy=0,\quad i=1,2,\quad j=1,2,3,4,5,$$ we obtain
\begin{equation*}
|\xi_{1}|=O(\lambda_1^{\frac{3}{2}}),
\end{equation*}
and
\begin{equation*}
|\xi_{2}|=q+O(\lambda_2^{\frac{3}{2}}).
\end{equation*}
\subsection{The selection of $\lambda_i$ and $\xi_i$}\ \\
Recall the orthogonal condition
\begin{equation*}
\begin{split}
\int_{B_{2R}} H_iZ_j=0\ \ \ \ i=1,2\ \ \ \ j=1,2,3,4,5,6,\\
\end{split}
\end{equation*}
For given $\psi$, define $\lambda_i$, $\xi_i$ as the unique solution of the above orthogonal condition. We require
\begin{equation}
\lambda_i(T)=0,\ \ \ \ \xi_1(T)=0,\ \ \ \ \xi_2(T)=q.
\end{equation}
Divide $\lambda_i$ into
\begin{equation*}
\lambda_i=\lambda_{i,0}+\lambda_{i,1}.
\end{equation*}
For $\lambda_1$
\begin{equation}
\begin{split}
(\dot{\lambda}_{1,0}+\dot{\lambda}_{1,1})\int_{B_{2R}}W_{6}^2-
\lambda_1^\frac{1}{2}\int_{B_{2R}}pU(y_1)^{p-1}W_{6}(y_1)(Z_1+Z_2-\psi)dy_1=0.
\end{split}
\end{equation}
Then by the choosing of $\lambda_{1,0}$, we have
\begin{equation*}
\dot{\lambda}_{1,0}+\kappa_1(T-t)\lambda_{1,0}^{\frac{1}{2}}=0.
\end{equation*}
Substitute it into equation (3.6)
\begin{equation}
\begin{split}
&\dot{\lambda}_{1,1}\int_{R^5}W_{6}^2-\lambda_1^\frac{1}{2}\int_{R^5}pU^{p-1}W_{6}M(T-t)dy_1\\
&+
\lambda_{1,0}^\frac{1}{2}\int_{R^5}pU^{p-1}W_{6}M(T-t)dy_1-\int_{R^5/B_{2R}}Z_{6}^2(\dot{\lambda}_{1,0}
+\dot{\lambda}_{1,1})dy_1\\
&-
\lambda_1^\frac{1}{2}\int_{B_{2R}}pU(y_1)^{p-1}W_{6}(y_1)[M(T-t)(-\frac{|z|^2}{10})+Z_2-\psi]dy_1=0.
\end{split}
\end{equation}
That is
\begin{equation}
\begin{split}
&\dot{\lambda}_{1,1}\int_{R^5}W_{6}^2+\frac{M(T-t)\lambda_{1,1}}{(\lambda_1^\frac{1}{2}+
\lambda_{1,0}^\frac{1}{2})}\int_{R^5}pU^{p-1}W_{6}dy_1-\int_{R^5/B_{2R}}W_{6}^2(\dot{\lambda}_{1,0}
+\dot{\lambda}_{1,1})dy_1\\
&-\lambda_1^\frac{1}{2}\int_{B_{2R}}pU(y_1)^{p-1}W_{6}(y_1)[M(T-t)(-\frac{|z|^2}{10})+Z_2-\psi]dy_1=0.
\end{split}
\end{equation}
Because in $y_1\in B_{2R}$, $|z|\leq 2R(T-t)^{\frac{7}{2}}$, then
\begin{equation}
\lambda_1^\frac{1}{2}\int_{B_{2R}}pU(y_1)^{p-1}W_{6}(y_1)[M(T-t)(-\frac{|z|^2}{10})]\leq \lambda_1^\frac{1}{2} (T-t)^8.
\end{equation}
Since $Z_2$ is odd, thus
\begin{equation*}\begin{split}
&\int_{B_{2R}}pU(y_1)^{p-1}W_{6}(y_1)Z_2(x+\xi_1)dy_1=\int_{B_{2R}}pU(y_1)^{p-1}W_{6}(y_1)Z_2(y_1\lambda_1)dy_1=0.
\end{split}
\end{equation*}
Therefore
\begin{equation}
\begin{split}
&|\lambda_1^\frac{1}{2}\int_{B_{2R}}pU(y_1)^{p-1}W_{6}(y_1)Z_2dy_1|\\&=
|\lambda_1^\frac{1}{2}\int_{B_{2R}}pU(y_1)^{p-1}W_{6}(y_1)[Z_2(x)-Z_2(y_1\lambda_1)]dy_1|\\
&\lesssim|\xi_1\lambda_1^\frac{1}{2}\int_{B_{2R}}pU(y_1)^{p-1}W_{6}(y_1)\nabla Z_2(s)dy_1|\\
&\lesssim |\xi_1\lambda_1^\frac{1}{2}\int_{B_{2R}}pU(y_1)^{p-1}W_{6}(y_1)dy_1|.
\end{split}
\end{equation}
By (3.9)-(3.10) and the definition of $\psi$ in  (2.5), equation (3.8) can be written as
\begin{equation}
\begin{split}
&\dot{\lambda}_{1,1}\int_{R^5}W_{6}^2+\int_{R^5}U^{p}dy_1M(T-t)\frac{\lambda_{1,1}}{2\lambda_{1,0}^\frac{1}{2}}+\delta_{R,T,\delta_0}(\dot{\lambda}_{1,0}+\dot{\lambda}_{1,1}
+(T-t)\lambda_1^{\frac{1}{2}})\\
&-\int_{R^5}U(y_1)^{p}dy_1M(T-t)\frac{\lambda_{1,1}^2}{2\lambda_{1,0}^\frac{1}{2}(\lambda_1^\frac{1}{2}+\lambda_{1,0}^\frac{1}{2})^2}=0,
\end{split}
\end{equation}
where when $R\rightarrow\infty$, $T\rightarrow0$, $\delta_0\rightarrow0$, $\delta_{R,T,\delta_0}\rightarrow0$. To solve $\lambda_{1,1}$, we need the following lemma.
\begin{lemma}
For sufficient small $\epsilon>0$, consider
\begin{equation}
\begin{cases}
\dot{\lambda}+\epsilon \frac{\lambda}{T-t}=(T-t)^3 h,\\
\lambda(T)=0.
\end{cases}
\end{equation}
If $\|h\|_{0,\alpha}<\infty$, (3.13) has a solution $\lambda$, which can be write as
\begin{equation}
\lambda=-(T-t)^{\epsilon}\int_t^T (T-s)^{3-\epsilon}h(s)ds,
\end{equation}
and there exists constant $C_2$ satisfies
\begin{equation}
\|\lambda\|_{1,\alpha}\leq C_2\|h\|_{0,\alpha}.
\end{equation}
\end{lemma}
\noindent Proof: In fact, if $\epsilon$ is small enough, $\int_t^T (T-s)^{3-\epsilon}h(s)ds$ is integrable, then $$(T-t)^{\epsilon}\int_t^T (T-s)^{3-\epsilon}h(s)ds$$ is well defined. One can check directly that (3.13) satisfies (3.12).\\
Then we check (3.14). By equation (3.12), $$\dot{\lambda}=-\epsilon \frac{\lambda}{T-t}+(T-t)^3 h,$$ where $(T-t)^3 h$ is in $C^\alpha$. By (3.13) $$-\epsilon \frac{\lambda}{T-t}=\epsilon(T-t)^{\epsilon-1}\int_t^T (T-s)^{3-\epsilon}h(s)ds,$$ it's easy to check it is in $C^1$ and $$|\frac{d(T-t)^{\epsilon-1}\int_t^T (T-s)^{3-\epsilon}h(s)ds}{dt}|\leq \|h\|_{\infty},$$ then Lemma 3.1 follows  immediately.

After finished the proof of Lemma 3.1. We may rewrite (3.11) as a linear equation about $\overline{\lambda_{1,1}}$
\begin{equation}
\begin{split}
&\dot{\overline{\lambda}}_{1,1}\int_{R^5}W_{6}^2+\frac{2M\int_{R^5}U^{p}dy_1}{\kappa_1(T-t)}\overline{\lambda_{1,1}}+\delta_{R,T,\delta_0}(\dot{\lambda}_{1,0}+\dot{\lambda}_{1,1}
+(T-t)\lambda_1^{\frac{1}{2}})\\
&-\int_{R^5}U(y_1)^{p}dy_1M(T-t)\frac{\lambda_{1,1}^2}{2\lambda_{1,0}^\frac{1}{2}(\lambda_1^\frac{1}{2}+\lambda_{1,0}^\frac{1}{2})^2}=0.
\end{split}
\end{equation}
If M sufficient small, then for given $\lambda_{1,1}\in C^{1,\alpha}$ with $|\lambda_{1,1}|\leq \delta \lambda_{1,0}, |\dot{\lambda}_{1,1}|\leq \delta |\dot{\lambda}_{1,0}|$, $\delta$ small enough. By Lemma 3.1, we can solve $\overline{\lambda_{1,1}}$
with the following estimate
\begin{equation*}
\|\overline{\lambda_{1,1}}\|_{1,\alpha}\leq C_1(\delta_{R,T,\delta_0}+\|\lambda_{1,1}\|_{1,\alpha}^2)\lambda_{1,0}.
\end{equation*}
Due to our choice of R sufficiently large, T sufficiently small, $\delta_0$ sufficiently small and $\delta$ small, by the fixed point theorem we can solve $\lambda_1$.

Similarly, we can solve $\lambda_2$ in a same way. As for $\xi_i=(\xi_{i,1},\cdots,\xi_{i,5})$, and $\xi_{i,j}$, defined by following equation
\begin{equation*}
\int_{B_{2R}}H_i(y_i,t)W_j(y_i)dy_i=0.
\end{equation*}
That is, for $j=1,\cdots,5$,
\begin{equation*}
\int_{B_{2R}} \lambda_i\dot{\xi}_i(\partial_j U(y_i))^2dy_i+\int_{B_{2R}}\lambda_i^{\frac{3}{2}}pU(y_i)^{p-1}(-Z_1-Z_2+\psi(\xi_i+\lambda y_i))dy_i=0.
\end{equation*}
Because $Z_1$, $Z_2$ are given, then for fixed $\psi$, after solve $\lambda_i(\psi)$ we can obtain $\xi_{i,j}$ from the above equations.

\section{The inner problem}\label{chap:5}
\noindent In this section, we use Lemma 2.1 to obtain an estimate of $\psi_1$ and $\psi_2$. Recall the inner problem
\begin{equation*}\begin{cases}
\lambda_1^2\partial_t\phi_{1}=\Delta_{y_1}\phi_1+pU^{p-1}(y_1)+H_1(\psi,\lambda_1,\lambda_2,\xi_1,\xi_2)\ \ \ (y_1,t)\in\ B_{2R}(0)\times(0,T),\\
\phi_{1}(y,0)=\phi_{1,0} \ \ \ \ \ \ \ \ \ \ \ \ \ \ \ \ \ \ \ \ \ \ \ \ \ \ \ \ \ \ \ \ \ \ \ \ \ \ \ \ \ \ \ \ \ \ \ \ \ \ \ \ \ \ \ y_1\in B_{2R}.
\end{cases}
\end{equation*}
\begin{equation*}\begin{cases}
\lambda_2^2\partial_t\phi_{2}=\Delta_{y_2}\phi_2+pU^{p-1}(y_2)+H_2(\psi,\lambda_1,\lambda_2,\xi_1,\xi_2)\ \ \ (y_2,t)\in\ B_{2R}(0)\times(0,T),\\
\phi_{2}(y,0)=\phi_{2,0}\ \ \ \ \ \ \ \ \ \ \ \ \ \ \ \ \ \ \ \ \ \ \ \ \ \ \ \ \ \ \ \ \ \ \ \ \ \ \ \ \ \ \ \ \ \ \ \ \ \ \ \ \ \ \ y_2\in B_{2R}.
\end{cases}
\end{equation*}
Where
\begin{equation*}
H_i(\psi,\lambda_i,\xi_i):=\lambda_i^{\frac{3}{2}}pU(y_i)^{p-1}(-Z_1(\xi_i+\lambda_i y_i)-Z_2(\xi_i+\lambda_i y_i)+\psi(\xi_i+\lambda_i y_i))+E_i(y_i,t),
\end{equation*}

\noindent For given $\psi$ satisfies (2.5), then
\begin{equation*}
H_1(\psi)\leq\lambda_1^{\frac{7}{4}}\frac{1}{1+|y_1|^3}.
\end{equation*}
\begin{equation*}
H_2(\psi)\leq\lambda_2^{\frac{3}{2}}\frac{1}{1+|y_2|^3}.
\end{equation*}
Notice that for any $a\in(0,1)$, the decay rate of $H_i(y_i,t)$ is 2+a. Also due to the choice of $\lambda_i$ and $\xi_i$, $H_i$ satisfies the orthogonality condition. Then by Lemma 2.1 we obtain
\begin{proposition}For any $a\in(0,1)$\\
\begin{equation*}\|\phi_1(y_1,t)\|_{*a,\frac{7}{4}}\leq C_3\|H_1(y_1,t)\|_{2+a,\frac{7}{4}}.\end{equation*}
 \begin{equation*}\|\phi_2(y_2,t)\|_{*a,\frac{3}{2}}\leq C_3\|H_2(y_1,t)\|_{2+a,\frac{3}{2}}.\end{equation*}
\end{proposition}
\noindent The norm here means in the corresponding variable, that is
\begin{equation}
|\phi_1(y_1,t)|\leq \lambda_1^{\frac{7}{4}}\frac{R^{6-a}}{1+|y_1|^6},\ \ \ \ |\phi_2(y_2,t)|\leq \lambda_2^{\frac{3}{2}}\frac{R^{6-a}}{1+|y_2|^6}.
\end{equation}

\section{The outer problem}\label{chap:6}
\noindent To solve outer problem, we also need to treat it as a linear problem
\begin{equation}\begin{cases}
\partial_t \psi_{*}&=\Delta_x\psi_*+g(\psi)\ \ \ \ (x,t)\in\  \mathbb{R}^5\times (0,T),\\
\psi_*(\cdot,0)&=\psi_0\ \ \ \ \ \ \ \ \ \ \ \ \ \ \ \ \ \ \ \ x\in\  \mathbb{R}^5
\end{cases}
\end{equation}
Our goal of this section is to prove $\psi$ has the following estimates
\begin{equation*}
\begin{split}
|\psi|\leq\begin{cases}\delta_0 (T-t)(1+|z|^4)\ &\ |z|\leq(T-t)^{-\frac{1}{4}},\\
\frac{\delta_0}{1+|x|^2}\ &\ |z|>(T-t)^{-\frac{1}{4}},
\end{cases}
\end{split}
\end{equation*}
where
\begin{equation*}\delta_0=\frac{\|Z_2\|_\infty}{10}.\end{equation*}
Define space $B$ with the norm
\begin{equation*}
\|\psi\|_B:=\sup_{x\in R^5,t\in[0,T)}|\psi(x,t)|(T-t)^{-1}(1+|z|^4),\ \ \ \ z=\frac{x}{T-t},  \ \ \ \  x\in \mathbb{R}^5.
\end{equation*}
From classic parabolic theorem, for given $\psi$ satisfies $\|\psi\|_B<\infty$, equation (5.1) has a unique solution $\psi_{*}$. Define
\begin{equation*}
\psi_*=T[\psi].
\end{equation*}

\subsection{Choosing of the parameters}\ \\
We choose initial data of this form
\begin{equation}\psi_0=(\boldsymbol{d}\cdot \boldsymbol{e})\eta(\frac{x}{e^{\frac{3}{8}\tau}}),\end{equation}
where $\boldsymbol{d}=(d_{\alpha_0},d_{\alpha_1}\cdots,d_{\alpha_m})$ is a parameter and $\boldsymbol{e}=(e_{\alpha_0},e_{\alpha_1}\cdots,e_{\alpha_m})$  is chosen such that $\boldsymbol{e}=(e_{\alpha_0},e_{\alpha_1}\cdots,e_{\alpha_m})$ are the eigenfunctions of $-A_z$ with eigenvalue
less than or equal 2.\\
We use the self-similar variable
\begin{equation*}
\varphi(z,\tau)=\psi_*(x,t),\ \ \ \ T-t=e^{-\tau}.\\
\end{equation*}
The function $\varphi(z,\tau)$ solves
\begin{equation}
\begin{split}
\begin{cases}
\varphi_\tau=\textit{A}\varphi+e^{-\tau}G(\phi,\psi,\lambda_1,\lambda_2,\xi_1,\xi_2)\ \ \ \ \ \ \ \ \ (z,\tau)\in\  \mathbb{R}^5\times (\tau_0,\infty)\\
\varphi\mid_{\tau=\tau_0}=(\boldsymbol{d}\cdot \boldsymbol{e})\eta(\frac{z}{e^{\frac{3}{8}\tau}}).
\end{cases}
\end{split}
\end{equation}
We decompose the initial data in the eigenspace of $A$,  into two parts, $Y_m=span\{e_{\alpha_0},e_{\alpha_1}\cdots,e_{\alpha_m}\}$ and it's orthogonal complement,
\begin{equation*}
(\boldsymbol{d}\cdot \boldsymbol{e})\eta(\frac{z}{e^{\frac{3}{8}\tau}})=\sum d_j(e_j\eta(\frac{z}{e^{\frac{3}{8}\tau}}),e_k)e_k+\{(\boldsymbol{d}\cdot \boldsymbol{e})\eta(\frac{z}{e^{\frac{3}{8}\tau}})\}^\perp,
\end{equation*}
especially
\begin{equation*}
(\{(\boldsymbol{d}\cdot \boldsymbol{e})\eta(\frac{z}{e^{\frac{3}{8}\tau}})\}^\perp,e_i)=0,\ \ \ \ i=1,\cdots,m,\ \ \ \tau=\tau_0.
\end{equation*}
We define $\Phi$ as
\begin{equation*}
\varphi=\boldsymbol{b}(\tau)\cdot \boldsymbol{e}+\Phi.
\end{equation*}
where $\boldsymbol{b}(\tau)=(b_{\alpha_1}(\tau),\cdots,b_{\alpha_m}(\tau))$.
Choose $b_{\tau}$ as
\begin{equation}
b_{\alpha_i}(\tau)=-e^{-\lambda_{\alpha_i} \tau}\int_{\tau}^{\infty} e^{(\lambda_{\alpha_i}-1)\tau'}(G,e_{\alpha_i})_\rho d\tau'.
\end{equation}
Then $\Phi$ satisfies
\begin{equation}
\Phi_\tau=\textit{A}_z \Phi+e^{-\tau}G-\sum_{i=1}^m\lambda_{\alpha_i} b_{\alpha_i} e_{\alpha_i}-\sum_{i=1}^m\frac{db_{\alpha_i}}{d\tau}e_{\alpha_i}:=\textit{A}_z \Phi+e^{-\tau}G^{\perp}.
\end{equation}
By Lemma 5.1, it's easy to see
\begin{equation}
|(G,e_{\alpha_i})|\lesssim e^{-\frac{7}{6}\tau}.
\end{equation}
From the definition of $\boldsymbol{d}$ we have
\begin{equation}
(I+D)\boldsymbol{d}=\boldsymbol{b}(\tau_0),
\end{equation}
where $D$ is a $(m+1)\times(m+1)$ matrix
\begin{equation}
D_{ij}=(e_{\alpha_i},(1-\eta(\frac{z}{e^{\frac{3}{8}\tau}})\mid_{\tau=\tau_0}e_{\alpha_j}).
\end{equation}
For $\tau_0$ sufficient large, $|D_{ij}|<<1$, then
\begin{equation}
|\boldsymbol{d}|\lesssim |\boldsymbol{b}(\tau)|\lesssim e^{-\frac{13}{6}\tau_0}.
\end{equation}
Because $e_\alpha$ are given functions with eigenvalues less than or equal to 2, so
\begin{equation}
\begin{split}
|b_\alpha e_\alpha|\lesssim e^{-\frac{13}{6}\tau}(1+|z|^4).
\end{split}
\end{equation}
\subsection{The estimate of G}\ \\
We give an estimate of G in this subsection. Let $R=e^{\frac{\tau_0}{2}}$. Recall (2.11) $G$ could be written as
\begin{equation*}
\begin{split}
G=&\sum_{i=1,2} p\lambda_i^{-2}(1-\eta_{R,i})U(y_i)^{p-1}(-Z_1\eta-Z_2+\psi)+\sum_{i=1,2} A[\phi_i]\\
&+\sum_{i=1,2} B[\phi_i]+\sum_{i=1,2} \lambda_i^{-\frac{7}{2}}E_i(1-\eta_{R,i})+N.
\end{split}
\end{equation*}
Let $\textbf{1}_{z\in \Omega}(z)$ be a function on $\mathbb{R}^5$ defined by  $\textbf{1}_{z\in \Omega}(z)=1$ if $z\in\Omega$ and  $\textbf{1}_{z\in \Omega}(z)=0$ if $z\notin \Omega$.
\begin{lemma} Under the above notation, we have
\begin{equation*}
\begin{split}
&|G|\\ &\lesssim
\frac{\lambda_1^{-2}}{1+|y_1|^{2+a}}e^{-\tau}R^{-\frac{1}{2}}\textbf{1}_{{|x|}\leq 1}+\frac{\lambda_2^{-2}}{1+|y_2|^{2+a}}R^{-\frac{1}{2}}\textbf{1}_{{|x-q|}\leq 1}+e^{-\frac{3}{4}\tau}|z|^4\textbf{1}_{e^{\frac{3}{8}\tau}\leq{|z|}\leq 2e^{\frac{3}{8}\tau}}\\&+e^{-\frac{7\tau}{6}}|z|^{4}\textbf{1}_{{|z|}\leq 2e^{\frac{\tau}{2}}}
+e^{-\frac{7\tau}{3}}|z|^{10}\textbf{1}_{{|z|}\leq 2e^{\frac{\tau}{4}}}+e^{-\frac{7\tau}{3}}|z|^{6}\textbf{1}_{{|z|}\leq 2e^{\frac{3}{8}\tau}}+e^{-\frac{7}{3}\tau}\textbf{1}_{{|z|}\leq 1}\\
&+\frac{\delta_0^{\frac{7}{3}}}{1+|x|^{4}}\textbf{1}_{{|z|}\geq e^{\frac{\tau}{4}}}+
\frac{\|Z_2\|_\infty^p}{1+|x|^4}\textbf{1}_{{|z|}\geq e^{\frac{\tau}{2}}}+\frac{e^{-2\tau}}{1+|x|^{3}}\textbf{1}_{{|x|}\geq 1}\\
&\lesssim\lambda_1^{-2}\frac{1}{1+|y_1|^{2+a}}e^{-\tau}R^{-\frac{1}{2}}\textbf{1}_{{|x|}\leq 1}+e^{-5\tau}|z|^{18}+e^{-\frac{3}{4}\tau}|z|^4\textbf{1}_{e^{\frac{3}{8}\tau}\leq{|z|}\leq 2e^{\frac{3}{8}\tau}}\\&+e^{-\frac{7\tau}{6}}|z|^{4}\textbf{1}_{{|z|}\leq 2e^{\frac{\tau}{2}}}
+e^{-\frac{7\tau}{3}}|z|^{10}\textbf{1}_{{|z|}\leq 2e^{\frac{\tau}{4}}}+e^{-\frac{7\tau}{3}}|z|^{6}\textbf{1}_{{|z|}\leq 2e^{\frac{3}{8}\tau}}\\
&+e^{-\frac{7}{3}\tau}\textbf{1}_{{|z|}\leq 1}
+\frac{\delta_0^{\frac{7}{3}}}{1+|x|^{4}}\textbf{1}_{{|z|}\geq e^{\frac{\tau}{4}}}+
\frac{\|Z_2\|_\infty^p}{1+|x|^4}\textbf{1}_{{|z|}\geq e^{\frac{\tau}{2}}}+\frac{e^{-2\tau}}{1+|x|^{3}}\textbf{1}_{{|x|}\geq 1}.
\end{split}
\end{equation*}
\end{lemma}
\noindent Proof: We can estimate the terms of G one by one.\\
For $|z|\leq e^{\frac{\tau}{2}}$, then
\begin{equation*}
\begin{split}
&\lambda_1^{-2}(1-\eta_{R,1})U(y_1)^{p-1}(-Z_1\eta-Z_2+\psi)\\
&\lesssim\lambda_1^{-2}(1-\eta_{R,1})\frac{1}{y_1^4}(e^{-\tau}
+e^{-\frac{\tau}{2}}|z|)\textbf{1}_{{|z|}\leq 1}+\lambda_1^2(1-\eta_{R,1})\frac{1}{y_1^3}e^{-\frac{7\tau}{2}}\textbf{1}_{1\leq|z|\leq e^{\frac{\tau}{2}}}\\
&\lesssim\lambda_1^{-2}\frac{1}{y_1^4}(1-\eta_{R,1})e^{-\tau}\textbf{1}_{{|z|}\leq1}+\lambda_1^{-2}\frac{1}{y_1^3}(1-\eta_{R,1})e^{-\frac{\tau}{2}}e^{-\frac{7\tau}{2}}\textbf{1}_{{|z|}\leq 1}+\lambda_1^2(1-\eta_{R,1})\frac{1}{y_1^3}e^{-\frac{7\tau}{2}}\textbf{1}_{1\leq|z|\leq e^{\frac{\tau}{2}}}\\
&\lesssim\lambda_1^{-2}\frac{1}{1+|y_1|^{2+a}}e^{-\tau}R^{-1}\textbf{1}_{{|z|}\leq 1}+\lambda_1^{-2}\frac{1}{1+|y_1|^{2+a}}e^{-\tau}R^{-1}\textbf{1}_{{|x|}\leq 1}\textbf{1}_{1\leq|z|\leq e^{\frac{\tau}{2}}}.
\end{split}
\end{equation*}
If $|z|\geq e^{\frac{\tau}{2}}$, then
\begin{equation*}
\begin{split}
\lambda_1^{-2}(1-\eta_{R,1})U(y_1)^{p-1}(-Z_1\eta_1-Z_2+\psi)&\lesssim\frac{e^{-2\tau}}{1+|x|^{3}}\textbf{1}_{{|x|}\geq 1}.
\end{split}
\end{equation*}
Because the control of $Z_1$, $Z_2$ and $\psi$, we deduce
\begin{equation*}
\begin{split}
&\lambda_2^{-2}(1-\eta_{R,2})U(y_2)^{p-1}(-Z_1\eta_1-Z_2+\psi)\\
&\lesssim\lambda_2^2(1-\eta_{R,2})\frac{1}{y_2^4}(e^{-\tau}+e^{-\frac{\tau}{2}}|z|)\textbf{1}_{{|x-q|}\leq 1}+\frac{e^{-2\tau}}{1+|x|^{3}}\textbf{1}_{{|x-q|}\geq 1}\\
&\lesssim\lambda_2^{-2}(1-\eta_{R,2})\frac{1}{y_2^4}\textbf{1}_{{|x-q|}\leq 1}+\frac{e^{-2\tau}}{1+|x|^{3}}\textbf{1}_{{|x-q|}\geq 1}\\
&\lesssim\lambda_2^{-2}\frac{1}{1+|y_2|^{2+a}}R^{-1}\textbf{1}_{{|x-q|}\leq 1}+\frac{e^{-2\tau}}{1+|x|^{3}}\textbf{1}_{{|x-q|}\geq 1}.\\
\end{split}
\end{equation*}
In the self-similar variable
\begin{equation*}
\lambda_2^{-2}(1-\eta_{R,2})\frac{1}{y_2^4}\textbf{1}_{{|x-q|}\leq 1}\lesssim e^{-4\tau}+e^{-5\tau}|z|^{18}.
\end{equation*}
From Proposition 4.1
\begin{equation*}
|\phi_1(y_1,t)|\leq C_1 \lambda_1^{\frac{7}{4}}\frac{R^{6-a}}{1+|y_1|^6},\quad |\phi_2(y_2,t)|\leq C_2 \lambda_2^{\frac{3}{2}}\frac{R^{6-a}}{1+|y_2|^6},
\end{equation*}
therefore
\begin{equation*}
\begin{split}
A[\phi_1]+B[\phi_1]&=\lambda_1^{-\frac{7}{2}}[\Delta_{y_1}\eta_{R,1}\phi_1+2\nabla_{y_1}\eta_{R,1}\nabla_{y_1}\phi_1]\\
&\lesssim\lambda_1^{-2}e^{-\tau}R^{-2-a}\textbf{1}_{R\leq{|y_1|}\leq 2R}+\frac{e^{-8\tau}R^{6-a}}{\lambda_1^2(1+|y_1|^6)}\textbf{1}_{{|y_1|}\leq 2R}\\
&\lesssim\lambda_1^{-2}\frac{1}{1+|y_1|^{2+a-\frac{1}{2}}}e^{-\tau}R^{-\frac{1}{2}}\textbf{1}_{R\leq{|y_1|}\leq 2R}+\lambda_1^{-2}\frac{1}{1+|y_1|^{2+a}}e^{-\tau}R^{-1}\textbf{1}_{{|y_1|}\leq 2R}.
\end{split}
\end{equation*}
Similarly
\begin{equation*}
\begin{split}
A[\phi_2]+B[\phi_2]&\lesssim\lambda_2^{-2}R^{-2-a}\textbf{1}_{R\leq{|y_2|}\leq 2R}+\lambda_1^{-2}\frac{1}{1+|y_1|^{2+a}}e^{-\tau}R^{-\frac{1}{2}}\textbf{1}_{{|y_2|}\leq 2R}\\
&\lesssim\lambda_2^{-2}\frac{1}{1+|y_2|^{2+a-\frac{1}{2}}}R^{-\frac{1}{2}}\textbf{1}_{R\leq{|y_2|}\leq 2R}+\lambda_1^{-2}\frac{1}{1+|y_1|^{2+a}}e^{-\tau}R^{-\frac{1}{2}}\textbf{1}_{{|y_2|}\leq 2R}\\
&\lesssim e^{-4\tau}+e^{-5\tau}|z|^{18}.
\end{split}
\end{equation*}
Use the definition of $E_1$, we have
\begin{equation*}
\begin{split}
\lambda_1^{-\frac{7}{2}}E_1(1-\eta_{R,1})
&\lesssim\lambda_1^{-2}\lambda_1^{\frac{1}{2}}e^\tau\frac{1}{|y_1|^3}(1-\eta_{R,1})\textbf{1}_{|x|\leq1}+\lambda_1^{-2}\lambda_1^{\frac{1}{2}}e^\tau\frac{1}{|y_1|^3}(1-\eta_{R,1})\textbf{1}_{|x|\geq1}\\
&\lesssim\lambda_1^{-2}\frac{1}{1+|y_1|^{2+a}}e^{-\tau}R^{-1}\textbf{1}_{|x|\leq1}+\frac{e^{-5\tau}}{1+|x|^{3}}\textbf{1}_{|x|\geq1}.
\end{split}
\end{equation*}
Similarly
\begin{equation*}
\begin{split}
\lambda_2^{-\frac{7}{2}}E_2(1-\eta_{R,2})
&\leq\frac{\lambda_2^{-2}\lambda_2^{\frac{1}{2}}e^\tau }{|y_2|^3}(1-\eta_{R,2})\textbf{1}_{|x-q|\leq1}+\lambda_2^{-2}\lambda_2^{\frac{1}{2}}e^\tau\frac{1}{|y_2|^3}(1-\eta_{R,2})\textbf{1}_{|x-q|\geq1}\\
&\leq\frac{\lambda_2^{-2}}{1+|y_2|^{2+a}}R^{-1}\textbf{1}_{|x|\geq1}+\frac{e^{-2\tau}}{1+|x|^{3}}\textbf{1}_{|x-q|\geq1}\\
&\leq e^{-4\tau}+e^{-5\tau}|z|^{18}+\frac{e^{-2\tau}}{1+|x|^{3}}.
\end{split}
\end{equation*}
From the definition of $Z_1$ and $\eta_1$, it's easy to estimate the last part of $N$
\begin{equation*}
\begin{split}
|\frac{\partial \eta_1}{\partial t}Z_1+\nabla_x Z_1\nabla \eta_1-\Delta \eta_1 Z_1|
\lesssim e^{\frac{3}{4}\tau}\textbf{1}_{e^{\frac{3}{8}\tau}\leq{|z|}\leq 2e^{\frac{3}{8}\tau}}
\lesssim e^{-\frac{3}{4}\tau}|z|^4\textbf{1}_{e^{\frac{3}{8}\tau}\leq{|z|}\leq 2e^{\frac{3}{8}\tau}}.
\end{split}
\end{equation*}
Then we deal with another part of $N$. If $x$ is away from q,
\begin{equation*}
\begin{split}
&|U_{\lambda_1,\xi_1}+U_{\lambda_2,\xi_2}+\theta-Z_1 \eta_1-Z_2|^{p-1}(U_{\lambda_1,\xi_1}+U_{\lambda_2,\xi_2}+\theta-Z_1 \eta_1-Z_2)\\
&\quad -U_{\lambda_1,\xi_1}^p-U_{\lambda_2,\xi_2}^p-pU_{\lambda_1,\xi_1}^{p-1}(\theta-Z_1 \eta_1-Z_2)-pU_{\lambda_2,\xi_2}^{p-1}(\theta-Z_1 \eta_1-Z_2)\\
&\lesssim U_{\lambda_1\xi_1}^{p-2}(U_{\lambda_2,\xi_2}+\theta-Z_1 \eta_1-Z_2)^2+U_{\lambda_2\xi_2}^p+(U_{\lambda_2,\xi_2}+\theta-Z_1 \eta_1-Z_2)^p\\
&\quad +U_{\lambda_1\xi_1}^{p-1}U_{\lambda_2\xi_2}.
\end{split}
\end{equation*}
Because $\|\phi_1\|_{*a,\frac{7}{4}}$ is bounded
\begin{equation*}
\begin{split}
U_{\lambda_1\xi_1}^{p-2}(\lambda_1^{-\frac{3}{2}}\phi_1(y_1,t)\eta_R(y_1))^2
&\lesssim\lambda_1^{-\frac{1}{2}}\frac{1}{1+|y_1|}\frac{e^{-2\tau }R^{12-2a}}{1+|y_1|^{12}}\\
&\lesssim \lambda_1^{-2}\frac{1}{1+|y_1|^{13}}e^{-8\tau}R^{12-2a}\\
&\lesssim\lambda_1^{-2}\frac{1}{1+|y_1|^{2+a}}e^{-\tau}R^{-1}\textbf{1}_{{|x|}\leq 1}+\frac{e^{-2\tau}}{1+|x|^{3}}\textbf{1}_{{|x|}\geq 1}.
\end{split}
\end{equation*}
And we have
\begin{equation*}
\begin{split}
U_{\lambda_1\xi_1}^{p-2}Z_1^2 \eta_1^2
&\lesssim\lambda_1^{-2}\frac{\lambda_1^{\frac{3}{2}}e^{-2\tau}(1+|y_1|^{1+a})}{1+|y_1|^{2+a}}\textbf{1}_{{|z|}\leq 1}+e^{-2\tau}\eta_1\textbf{1}_{{|z|}\geq 1}
\\
& \lesssim\lambda_1^{-2}\frac{e^{-\frac{5}{2}\tau}}{1+|y_1|^{2+a}}\textbf{1}_{{|z|}\leq 1}+e^{-2\tau}\eta_1\textbf{1}_{{|z|}\geq 1}
\\
&\lesssim\lambda_1^{-2}\frac{1}{1+|y_1|^{2+a}}e^{-\tau}R^{-1}\textbf{1}_{{|z|}\leq 1}+e^{-2\tau}\eta_1\textbf{1}_{{|z|}\geq 1}.
\end{split}
\end{equation*}
\begin{equation*}
\begin{split}
U_{\lambda_1\xi_1}^{p-2}\Psi^2
\lesssim \lambda_1^{-2}\frac{\lambda_1^{\frac{3}{2}}e^{-2\tau}(1+|y_1|^{1+a})}{1+|y_1|^{2+a}}\textbf{1}_{|z|\leq e^{\frac{\tau}{4}}}+\frac{e^{-\frac{19}{8}\tau}}{1+|x|^4}\textbf{1}_{|z|\geq e^{\frac{\tau}{4}}}\\
\lesssim\lambda_1^{-2}\frac{1}{1+|y_1|^{2+a}}e^{-\tau}R^{-1}\textbf{1}_{|z|\leq e^{\frac{\tau}{4}}}+\frac{e^{-\frac{19}{8}\tau}}{1+|x|^4}\textbf{1}_{|z|\geq e^{\frac{\tau}{4}}}.
\end{split}
\end{equation*}
Similarly for $|z|\leq 1$,
\begin{equation*}
\begin{split}
U_{\lambda_1\xi_1}^{p-2}Z_2^2\lesssim\lambda^{-\frac{1}{2}}\frac{1}{1+|y_1|}e^{-\tau}&\lesssim\lambda_1^{-2}\frac{e^{-\frac{3}{2}\tau}}{1+|y_1|^{2+a}}\lesssim\lambda_1^{-2}\frac{1}{1+|y_1|^{2+a}}e^{-\tau}R^{-1}.
\end{split}
\end{equation*}
For $|z|\geq 1$, and since $Z_2$ is a fast decay function,
\begin{equation*}
\begin{split}
U_{\lambda_1\xi_1}^{p-2}Z_2^2\lesssim \frac{e^{-\frac{3}{2}\tau}}{1+|x|^4}.
\end{split}
\end{equation*}
And $U_{\lambda_1\xi_1}^{p-2}U_{\lambda_2\xi_2}^2$ can be control by $U_{\lambda_1\xi_1}^{p-1}U_{\lambda_2\xi_2}$ and $U_{\lambda_2\xi_2}^{p}$.

Moreover, for $|x|\leq 1$
\begin{equation*}
\begin{split}
U_{\lambda_1\xi_1}^{p-1}U_{\lambda_2\xi_2}\lesssim\lambda_1^{-2}\frac{e^{-3\tau}}{1+|y_1|^4}\textbf{1}_{|x|\leq1}.
\end{split}
\end{equation*}
For $|x|\geq 1$, $|x-q|\geq 1$, we have
\begin{equation*}
\begin{split}
U_{\lambda_1\xi_1}^{p-1}U_{\lambda_2\xi_2}\lesssim \frac{e^{-17\tau}}{1+|x|^4}.
\end{split}
\end{equation*}
\begin{equation*}
\begin{split}
(\lambda_1^{-\frac{3}{2}}\phi_1(y_1,t)\eta_R(y_1))^p&\lesssim\lambda_1^{-2}\frac{e^{-8\tau}e^{-\frac{7}{3}\tau}R^{12}}{1+|y_1|^{12}}\lesssim\frac{e^{-2\tau}}{1+|x|^{3}}\textbf{1}_{{|x|}\geq 1}.
\end{split}
\end{equation*}
\begin{equation*}
\begin{split}
(\lambda_2^{-\frac{3}{2}}\phi_2(y_2,t)\eta_R(y_2))^p\lesssim\lambda_2^{-2}\frac{e^{-4\tau}R^{12}}{1+|y_2|^{12}}.
\end{split}
\end{equation*}
As before
\begin{equation*}
\begin{split}
|\Psi|^p+|Z_1\eta_1|^p&\lesssim e^{-\frac{7\tau}{3}}(1+|z|^{\frac{28}{3}})\textbf{1}_{|z|\leq e^{\frac{\tau}{4}}}+\frac{\delta_0^{\frac{7}{3}}}{1+|x|^{\frac{14}{3}}}\textbf{1}_{|z|\geq e^{\frac{\tau}{4}}}+e^{-\frac{7\tau}{3}}(1+|z|^{\frac{14}{3}})\eta_1\\
&\lesssim e^{-\frac{7\tau}{3}}(1+|z|^{10})\textbf{1}_{|z|\leq e^{\frac{\tau}{4}}}+\frac{\delta_0^{\frac{7}{3}}}{1+|x|^{4}}\textbf{1}_{|z|\geq e^{\frac{\tau}{4}}}+e^{-\frac{7\tau}{3}}(1+|z|^{6})\eta_1.
\end{split}
\end{equation*}
Moreover
\begin{equation*}
\begin{split}
|Z_2|^p&\leq e^{-\frac{7\tau}{6}}|z|^{\frac{7}{3}}\textbf{1}_{|z|\leq e^{\frac{\tau}{2}}}+\frac{\|Z_2\|_\infty^p}{1+|x|^4}\textbf{1}_{|z|\geq e^{\frac{\tau}{2}}}\leq e^{-\frac{7\tau}{6}}|z|^{4}\textbf{1}_{|z|\leq e^{\frac{\tau}{2}}}+\frac{\|Z_2\|_\infty^p}{1+|x|^4}\textbf{1}_{|z|\geq e^{\frac{\tau}{2}}}.
\end{split}
\end{equation*}
If $|x-q|\leq 1$
\begin{equation*}
\begin{split}
&|U_{\lambda_1,\xi_1}+U_{\lambda_2,\xi_2}+\theta-Z_1 \eta_1-Z_2|^{p-1}(U_{\lambda_1,\xi_1}+U_{\lambda_2,\xi_2}+\theta-Z_1 \eta_1-Z_2)\\
&\quad -U_{\lambda_1,\xi_1}^p-U_{\lambda_2,\xi_2}^p-pU_{\lambda_1,\xi_1}^{p-1}(\theta-Z_1 \eta_1-Z_2)-pU_{\lambda_2,\xi_2}^{p-1}(\theta-Z_1 \eta_1-Z_2)\\
&\leq U_{\lambda_2,\xi_2}^{p-2}(U_{\lambda_1,\xi_1}+\theta-Z_1 \eta_1-Z_2)^2+U_{\lambda_1,\xi_1}^p+(U_{\lambda_1,\xi_1}+\theta-Z_1 \eta_1-Z_2)^p\\
&\quad +U_{\lambda_2,\xi_2}^{p-1}U_{\lambda_1,\xi_1}.
\end{split}
\end{equation*}
And
\begin{equation*}
\begin{split}
U_{\lambda_2,\xi_2}^{p-2}(U_{\lambda_1,\xi_1}+\theta-Z_1 \eta_1-Z_2)^2
\lesssim\frac{1}{\lambda_2^\frac{1}{2}(1+|y_2|)}
\lesssim\frac{e^{-\frac{\tau}{2}}}{\lambda_2^2(1+|y_2|^{2+a})}
\lesssim\frac{1}{\lambda_2^{2}(1+|y_2|^{2+a})}R^{-\frac{1}{2}}.
\end{split}
\end{equation*}
Also
\begin{equation*}
\begin{split}
U_{\lambda_1,\xi_1}^p+U_{\lambda_2,\xi_2}^{p-1}U_{\lambda_1,\xi_1}\lesssim e^{-14\tau}+\frac{e^{-6\tau}}{\lambda_2^2(1+|y|_2^{4})}\lesssim e^{-14\tau}+\lambda_2^{-2}\frac{1}{1+|y_2|^{2+a}}R^{-\frac{1}{2}}.
\end{split}
\end{equation*}
In conclusion, we complete the proof.

By using Lemma 5.1, we can easily deduce that
\begin{proposition} The function $G$ given by (2.11) satisfies
\begin{equation}
\|G\|_\rho\lesssim e^{-\frac{7}{6}\tau}.
\end{equation}
\end{proposition}
\noindent Proof:
Recall $y_1=\frac{x-\xi_1}{\lambda_1}$, $z=\frac{x}{\sqrt{T-t}}$ , that is $$z=\frac{(y_1\lambda_1+\xi_1)}{\sqrt{T-t}},$$ then we have
\begin{equation*}
\begin{split}
&|\int_{R^5}(\lambda_1^{-2}\frac{1}{1+|y_1|^{2+a}}e^{-\tau}R^{-\frac{1}{2}}\textbf{1}_{|x|\leq 1})^2 e^{-\frac{|z|^2}{4}} dz|
\\&=\int_{R^5}(\lambda_1^{-2}\frac{1}{1+|y_1|^{2+a}}e^{-\tau}R^{-\frac{1}{2}})^2 e^{-\frac{|z|^2}{4}}(\frac{\lambda_1}{\sqrt{T-t}})^5 dy_1\\
&\lesssim \lambda_1 e^{-2\tau}e^{-\frac{5}{2}\tau}\\
&\lesssim e^{-\frac{7}{2}\tau}.
\end{split}
\end{equation*}
And
\begin{equation*}
|\int_{R^5}(e^{-5\tau}|z|^{18})^2e^{-\frac{|z|^2}{4}} dz|\lesssim e^{-10\tau}.
\end{equation*}
For $\tau\geq\tau_0$ sufficient large, we have
\begin{equation*}
\begin{split}
|\int_{R^5}(e^{-\frac{3}{4}\tau}|z|^4\textbf{1}_{e^{\frac{3}{8}\tau}\leq{|z|}\lesssim 2e^{\frac{3}{8}\tau}})^2e^{-\frac{|z|^2}{4}} dz|&\lesssim e^{\frac{3}{2}\tau}e^{-\frac{e^{\frac{3}{4}\tau}}{4}}\leq e^{-5\tau}.
\end{split}
\end{equation*}
\begin{equation*}
|\int_{R^5}(e^{-\frac{7\tau}{6}}|z|^{4}\textbf{1}_{{|z|}\lesssim e^{\frac{\tau}{2}}})^2e^{-\frac{|z|^2}{4}} dz|\leq e^{-\frac{7\tau}{3}}.
\end{equation*}
\begin{equation*}
|\int_{R^5}(e^{-\frac{7\tau}{3}}|z|^{10}\textbf{1}_{{|z|}\lesssim e^{\frac{\tau}{4}}})^2e^{-\frac{|z|^2}{4}} dz|\leq e^{-\frac{14\tau}{3}}.
\end{equation*}
\begin{equation*}
|\int_{R^5}(e^{-\frac{7\tau}{3}}|z|^{6}\textbf{1}_{{|z|}\lesssim e^{\frac{3}{8}\tau}})^2e^{-\frac{|z|^2}{4}} dz|\leq e^{-\frac{14\tau}{3}}.
\end{equation*}
\begin{equation*}
|\int_{R^5}(e^{-\frac{7}{3}\tau}\textbf{1}_{{|z|}\leq 1})^2e^{-\frac{|z|^2}{4}} dz|\lesssim e^{-\frac{14\tau}{3}}.
\end{equation*}
Since the function $s^3e^{-\frac{s^2}{8}}$ is bounded and for $\tau\geq\tau_0$ sufficient large, we have
\begin{equation*}
\begin{split}
|\int_{R^5}(\frac{\delta_0^{\frac{7}{3}}}{1+|x|^{4}}\textbf{1}_{{|z|}\geq e^{\frac{\tau}{2}}})^2e^{-\frac{|z|^2}{4}} dz|
&\lesssim \int_{e^{\frac{\tau}{2}}}^\infty s^4 e^{-\frac{s^2}{4}}ds\\
&\lesssim \int_{e^{\frac{\tau}{2}}}^\infty s^3 e^{-\frac{s^2}{8}}se^{-\frac{s^2}{8}}ds\\
&\lesssim e^{-\frac{e^\tau}{8}}\\
&\lesssim e^{-5\tau}.
\end{split}
\end{equation*}
Similarly
\begin{equation*}
\begin{split}
|\int_{R^5}(\frac{e^{-2\tau}}{1+|x|^{3}}\textbf{1}_{{|x|}\geq 1}+\frac{\|Z_2\|_\infty^p}{1+|x|^4}\textbf{1}_{{|z|}\geq e^{\frac{\tau}{2}}})^2e^{-\frac{|z|^2}{4}} dz|\lesssim e^{-5\tau}.
\end{split}
\end{equation*}
Therefore
\begin{equation*}
\|G\|_\rho\lesssim e^{-\frac{7}{6}\tau}.
\end{equation*}
\subsection{The solving of $\Phi_1$}\ \\
We divide $\Phi$ into two parts, one part comes from the initial data, another comes from nonhomogeneous term, consider
\begin{equation}
\begin{split}
\begin{cases}
\partial_\tau\Phi_{1}=A_z\Phi_1\ \ \ \ \ \ \ \ \ \ \ \ \ \ \ \ \ \ \ \ \ \ \ \ \ (z,\tau)\in\mathbb{R}^5\times[\tau_0,\infty),\\
\Phi_1\mid_{\tau=\tau_0}=\{{(\boldsymbol{d}\cdot \boldsymbol{e})\eta(\frac{z}{e^{\frac{3}{8}\tau}})}\}^\perp.
\end{cases}
\end{split}
\end{equation}
\begin{lemma}
There exists constant $K_1$, such that
\begin{equation}
|\Phi_1(z,\tau)|\leq K_1e^{-2\tau}(1+|z|^4),\ \ \ \ \tau>\tau_0 +1.
\end{equation}
\end{lemma}
\noindent Proof: First we estimate the initial data
\begin{equation*}
\{{(\boldsymbol{d}\cdot \boldsymbol{e})\eta(\frac{z}{e^{\frac{3}{8}\tau}})}\}^\perp={(\boldsymbol{d}\cdot \boldsymbol{e})\eta(\frac{z}{e^{\frac{3}{8}\tau}})}
-\sum d_{\alpha_j}(e_{\alpha_j}\eta(\frac{z}{e^{\frac{3}{8}\tau}}),e_{\alpha_k})_\rho e_{\alpha_k},
\end{equation*}
As in the discuss of last section, we can deduce that
\begin{equation*}
|\{{(\boldsymbol{d}\cdot \boldsymbol{e})\eta(\frac{z}{e^{-\frac{3}{8}\tau}})}\}^\perp|\lesssim e^{-\frac{13}{6}\tau_0}+e^{-\frac{13}{6}\tau_0}|z|^4,
\end{equation*}
and
\begin{equation*}
\begin{split}
\|\{{(\boldsymbol{d}\cdot \boldsymbol{e})\eta(\frac{z}{e^{\frac{3}{8}\tau}})}\}^\perp\|_\rho\lesssim e^{-\frac{13}{6}\tau_0}.\\
\end{split}
\end{equation*}
We easily deduced that
\begin{equation*}
\begin{split}
\frac{1}{2}\partial_\tau\|\Phi_{1}\|_\rho^2=(A_z\Phi_1,\Phi_1)
\leq -\frac{5}{2}\|\Phi_{1}\|_\rho^2.
\end{split}
\end{equation*}
witch is
\begin{equation*}
\partial_\tau (\|\Phi_{1}\|_\rho^2 e^{5\tau})\leq 0.
\end{equation*}
Integrating it from $\tau_0$ to $\tau$
\begin{equation*}
\begin{split}
\|\Phi_1\|_\rho\leq\|{(\boldsymbol{d}\cdot \boldsymbol{e})\eta(\frac{z}{e^{\frac{3}{8}\tau}})}^\perp\|_\rho e^{-\frac{5}{2}(\tau-\tau_0)}\lesssim e^{-2\tau}.
\end{split}
\end{equation*}
Let $r_0$ be the constant such that in $|z|\geq r_0$, $e_2>0$. From the classical $L^2$ estimate of parabolic equation and the definition of $\rho$ we have, for $\tau>\tau_0 +1,|z|<r_0$
\begin{equation*}
|\Phi_1|\lesssim\|\Phi_1\|_{2,B_{2r_0}}\lesssim e^{-2\tau}.
\end{equation*}
Now we construct a comparison function $K_1e^{-2\tau}e_2(z)$. In fact, $K$ sufficient large, we prove it is a super solution.
When $|z|=r_0$
\begin{equation*}|\Phi_1|<K_1e^{-2\tau}e_2(z).\ \ \end{equation*}
For $\tau_0\leq\tau\leq\tau_0 +1$, by Lemma 2.3,
\begin{equation*}
|\Phi_1|\lesssim e^{-\frac{13}{6}\tau_0}+e^{-\frac{13}{6}\tau_0}e^{-2(\tau-\tau_0)}|z|^4.
\end{equation*}
Therefore when $\tau=\tau_0 +1$
\begin{equation*}
|\Phi_1|\lesssim e^{-\frac{13}{6}\tau_0}+e^{-\frac{13}{6}\tau_0}|z|^4\leq K_1e^{-2\tau}e_2(z).
\end{equation*}
We can easily verify
\begin{equation*}
\partial_\tau (Ke^{-2\tau}e_2(z))-A_z (Ke^{-2\tau}e_2(z))=0=\partial_\tau \Phi_1-A_z\Phi_1.
\end{equation*}
By comparison principle, $$|\Phi_1|<K_1e^{-2\tau}e_2(z).$$

\subsection{The solving of $\Phi_2$}\ \\
Recall the equation of $\Phi_2$
\begin{equation}
\begin{cases}
\partial_t\Phi_{2}-\Delta\Phi_2=G^\bot\ \ \ \ \ \ \ \ \ \ \ (z,\tau)\in\mathbb{R}^5\times[\tau_0,\infty),\\
\Phi_2\mid_{\tau=\tau_0}=0.
\end{cases}
\end{equation}
Our goal is to prove $\Phi_2$ vanishes around 0 and $\|\Phi_2\|_\infty$ small. So we first prove the following lemma.
\begin{lemma}
Let $A_z=\Delta_z-\frac{z}{2}\cdot\nabla$ be the operator as before, then for any $\tau\geq\tau_1\geq\tau$ we have
\begin{equation}
|\int_{\tau_1}^{\tau}e^{A_z(\tau-\tau')}e^{-\tau'}G^\bot|\lesssim \frac{e^{-\tau_1}}{R^{\frac{1}{2}}}+e^{-6\tau}|z|^{18}+e^{-\frac{7}{4}\tau}|z|^4+e^{-\frac{10}{3}\tau}|z|^{10}.
\end{equation}
\end{lemma}
\noindent Proof:
Notice
\begin{equation*}
\begin{split}
G^\bot=G-(G,e_{\alpha_i})e_{\alpha_i},
\end{split}
\end{equation*}
and
\begin{equation*}
\begin{split}
|b_\alpha e_\alpha|\lesssim e^{-\frac{13}{6}\tau}(1+|z|^4).
\end{split}
\end{equation*}
Then
\begin{equation*}
\begin{split}
G\lesssim&\lambda_1^{-2}\frac{1}{1+|y_1|^{2+a}}e^{-\tau}R^{-\frac{1}{2}}\textbf{1}_{{|x|}\leq 1}+e^{-5\tau}|z|^{18}+e^{-\frac{3}{4}\tau}|z|^4\textbf{1}_{e^{\frac{3}{8}\tau}\leq{|z|}\leq 2e^{\frac{3}{8}\tau}}\\
&+\frac{\delta_0^{\frac{7}{3}}}{1+|x|^{4}}\textbf{1}_{{|z|}\geq e^{\frac{\tau}{4}}}+\frac{\|Z_2\|_\infty^p}{1+|x|^4}\textbf{1}_{{|z|}\geq e^{\frac{\tau}{2}}}+\frac{e^{-2\tau}}{1+|x|^{3}}\textbf{1}_{{|x|}\geq 1}
+e^{-\frac{7\tau}{3}}(1+|z|^{10}).\\
\end{split}
\end{equation*}
We estimate the following terms
\begin{equation*}
\begin{split}
\frac{\delta_0^{\frac{7}{3}}}{1+|x|^{4}}\textbf{1}_{{|z|}\geq e^{\frac{\tau}{4}}}+\frac{\|Z_2\|_\infty^p}{1+|x|^4}\textbf{1}_{{|z|}\geq e^{\frac{\tau}{2}}}+\frac{e^{-2\tau}}{1+|x|^{3}}\textbf{1}_{{|x|}\geq 1}
\leq e^{-\tau}|z|^4.
\end{split}
\end{equation*}
Let $\Psi_1$ be the solution of following equation
\begin{equation*}
\begin{cases}
\partial_t\Psi_{1}-\Delta\Psi_1=\lambda_1^{-2}\frac{1}{1+|y_1|^{2+a}}e^{-\tau}R^{-\frac{1}{2}}\textbf{1}_{{|x|}\leq 1}\ \ \ \ \ \ \ \ \ \ \ (x,t)\in\mathbb{R}^5\times[0,T),\\
\Psi_1\mid_{t=0}=0.
\end{cases}
\end{equation*}
By Lemma 2.4, we obtain
\begin{equation*}
|\Psi_{1}|\lesssim \frac{e^{-\tau_1}}{R^{\frac{1}{2}}}.
\end{equation*}
$\Psi_{2}$ is the solution of following equation
\begin{equation*}
\begin{cases}
\partial_t\Psi_{2}-\Delta\Psi_1=&e^{-5\tau}|z|^{18}+e^{-\frac{3}{4}\tau}|z|^4\textbf{1}_{e^{\frac{3}{8}\tau}\leq{|z|}\leq 2e^{\frac{3}{8}\tau}}+\frac{\delta_0^{\frac{7}{3}}}{1+|x|^{4}}\textbf{1}_{{|z|}\geq e^{\frac{\tau}{4}}}\\
&+\frac{\|Z_2\|_\infty^p}{1+|x|^4}\textbf{1}_{{|z|}\geq e^{\frac{\tau}{2}}}
+e^{-\frac{7\tau}{3}}(1+|z|^{10})\ \ \ \ \ \ \ \ \ \ \ (x,t)\in\mathbb{R}^5\times[0,T),\\
\Psi_2\mid_{t=0}=0.
\end{cases}
\end{equation*}
By Duhamel's principle
\begin{equation*}
\begin{split}
|\Psi_{2}|\leq &\int_{\tau_1}^{\tau}e^{-\frac{10}{3}\tau'}(1+e^{-5(\tau-\tau')}|z|^{10})+e^{-6\tau'}(1+e^{-9(\tau-\tau')}|z|^{18})
\\
&+e^{-\frac{7}{4}\tau'}(1+e^{-2(\tau-\tau'}|z|^4)e^{-\frac{13}{6}\tau'}(1+e^{-2(\tau-\tau')}|z|^4)d\tau'\\
\lesssim &e^{-\frac{7}{4}\tau_1}+e^{-6\tau}|z|^{18}+e^{-\frac{7}{4}\tau}|z|^4+e^{-\frac{10}{3}\tau}|z|^{10}.
\end{split}
\end{equation*}
Together with $\Psi_{1}$,  we obtain the conclusion.

Next we improve the estimate of Lemma 5.4 when $\tau$ is large
\begin{lemma} Under the assumptions of Lemma 5.4, we have
\begin{equation}
\begin{split}
|\int_{\tau_0}^{\tau}e^{A_z(\tau-\tau')}e^{-\tau'}G^\bot d\tau'|&\lesssim \frac{e^{-\tau}}{R^\frac{1}{2}}+e^{-\frac{13}{6}\tau}e^{\frac{5}{4}\tau_0}
+e^{-\frac{7}{4}\tau}|z|^4+\frac{e^{-\tau}}{R^{\frac{1}{2}}}|z|^2
\\ &\quad +e^{-6\tau}|z|^{18}+e^{-\frac{10}{3}\tau}|z|^{10}.
\end{split}
\end{equation}
\end{lemma}
\noindent Proof:
We separate space-time into four cases:
\begin{equation}
\begin{split}
&1.\ z\in \mathbb{R}^5,\ \tau\in(\tau_0,\tau_0+1).\\
&2.\ |z|\leq 4,\ \tau\in(\tau_0+1,\infty).\\
&3.\ 4\leq|z|\leq e^{\frac{\tau-\tau_0}{2}},\ \tau\in(\tau_0+1,\infty).\\
&4.\ |z|>e^{\frac{\tau-\tau_0}{2}},\ \tau\in(\tau_0+1,\infty).
\end{split}
\end{equation}
\begin{itemize}
\item $z\in \mathbb{R}^5, \tau\in(\tau_0,\tau_0+1)$, we conclude by Lemma 5.4.\\
\item $|z|\leq 4$. Firstly we times $\Phi_2$ and integral with weight $\rho$ in (5.14)
\begin{equation*}
\frac{1}{2}\frac{d\|\Phi_2\|_{\rho}^2}{d\tau}=(A_z\Phi_2,\Phi_2)_{\rho}+e^{-\tau}(G^{\perp},\Phi_2)_{\rho}.
\end{equation*}
Since $(A_z\Phi_2,\Phi_2)_{\rho}\leq-\frac{5}{2}\|\Phi_2\|^2$, we deduce from Lemma 5.1 that
\begin{equation}
\|\Phi_2\|_{\rho}\lesssim e^{-\frac{13}{6}\tau}.
\end{equation}
Next, from Lemma 2.2 and the definition of $G^{\perp}$
\begin{equation*}
\begin{split}
|\int_{\tau_0}^{\tau-1}e^{A_z(\tau-\tau')}e^{-\tau'}G^{\bot} d\tau'|&=|\int_{\tau_0}^{\tau-1}e^{A_z\frac{1}{2}}e^{A_z(\tau-\tau'-\frac{1}{2})}e^{-\tau'}G^{\bot} d\tau'|\\
&\lesssim\int_{\tau_0}^{\tau-1}\frac{\frac{e^{-\frac{1}{2}|z|^2}}{e^{4(2-e^{-\tau_0}+e^{-\tau_0
-\frac{1}{2}})}}}{(e^{-\tau_0}-e^{-\tau_0-1})^{\frac{5}{4}}}e^{-\frac{5}{2}(\tau-\tau')}\|e^{-\tau'}G^{\bot} \|_\rho d\tau'\\
&\lesssim \int_{\tau_0}^{\tau-1}e^{\frac{5}{4}\tau_0}e^{-\frac{5}{2}(\tau-\tau')}e^{-\frac{13}{6}\tau'}d\tau'\\
&\lesssim e^{-\frac{13}{6}\tau}e^{\frac{5}{4}\tau_0}.\\
\end{split}
\end{equation*}
The integral from $\tau-1$ to $\tau$ can be estimated by Lemma 5.4, then
\begin{equation*}
\begin{split}
|\int_{\tau-1}^{\tau}e^{Az(\tau-\tau')}e^{-\tau'}G^{\bot}d\tau'|\lesssim\frac{e^{-\tau}}{R^\frac{1}{2}}.
\end{split}
\end{equation*}
\item $4\leq|z|\leq e^{\frac{\tau-\tau_0}{2}}, \tau >\tau_0+1$. Let $\tau_2$ define by $|z|=e^{\frac{\tau-\tau_2}{2}}$, from the definition, $\tau>\tau_2$,  a similarly calculation shows
\begin{equation*}
\begin{split}
|\int_{\tau_0}^{\tau_2}e^{Az(\tau-\tau')}e^{-\tau'}G^{\bot}d\tau'|&\leq |\int_{\tau_0}^{\tau_2}e^{\frac{5}{4}\tau_0}e^{-\frac{5}{2}(\tau-\tau')}e^{-\frac{13}{6}\tau'}G^{\bot}\tau'|\\
&\lesssim e^{\frac{1}{3}\tau_2}e^{-\frac{5}{2}\tau}e^{\frac{5}{4}\tau_0}\\
&\lesssim e^{-\frac{13}{6}\tau}e^{\frac{5}{4}\tau_0},\\
|\int_{\tau_2}^{\tau-1}e^{Az(\tau-\tau')}e^{-\tau'}G^{\bot}d\tau'|&\leq\frac{e^{-\tau_2}}{R^{\frac{1}{2}}}+e^{-6\tau}|z|^{18}+e^{-\frac{7}{4}\tau}|z|^4+e^{-\frac{10}{3}\tau}|z|^{10}\\
&\lesssim \frac{e^{-\tau}}{R^{\frac{1}{2}}}|z|^2+e^{-6\tau}|z|^{18}+e^{-\frac{7}{4}\tau}|z|^4+e^{-\frac{10}{3}\tau}|z|^{10}.\\
\end{split}
\end{equation*}
\item $|z|>e^{\frac{\tau-\tau_0}{2}}$. Similar as before, by Lemma 5.4 and $|z|>e^{\frac{\tau-\tau_0}{2}}$.
\begin{equation*}
\begin{split}
&|\int_{\tau_0}^{\tau}e^{Az(\tau-\tau')}e^{-\tau'}G^{\bot}d\tau'|\\
&\lesssim\frac{e^{-\tau_0}}{R^{\frac{1}{2}}}+e^{-6\tau}|z|^{18}+e^{-\frac{7}{4}\tau}|z|^4+e^{-\frac{10}{3}\tau}|z|^{10}\\
&\lesssim\frac{e^{-\tau}}{R^{\frac{1}{2}}}|z|^2+e^{-6\tau}|z|^{18}+e^{-\frac{7}{4}\tau}|z|^4+e^{-\frac{10}{3}\tau}|z|^{10}.
\end{split}
\end{equation*}
\end{itemize}
Combining all these estimates we complete the proof of Lemma 5.5.

Recall our goal is to show $\Phi_*(x,t)$ lays in space $X$ and it's norm is smaller than $\delta_0$. In the estimate of Lemma 5.5, when $\tau>\frac{15}{14}\tau_0$ the term $$e^{-\frac{13}{6}\tau}e^{\frac{5}{4}\tau_0}\in X,$$ and
\begin{equation*}
e^{-\frac{13}{6}\tau}e^{\frac{5}{4}\tau_0}< \delta_0 e^{-\tau}.
\end{equation*}
But when $\tau$ is near $\tau_0$ we can not obtain the same estimate. Thus we need to employ Lemma 5.4. In Lemma 5.4, let $\tau_1=\tau_0$, because $R=e^{\frac{\tau_0}{2}}$, the first term in Lemma 5.4 is
\begin{equation*}
\frac{e^{-\tau_0}}{R^{\frac{1}{2}}}={e^{-\tau_0}}e^{-\frac{\tau_0}{4}}={e^{-\frac{5}{4}\tau_0}}.
\end{equation*}
This is to say when $\tau<\frac{5}{4}\tau_0$ we have
\begin{equation}
|\Phi_2|<\delta_0 e^{-\tau}.
\end{equation}
So, in $|z|\leq e^{\frac{\tau}{4}}$, for any $\tau\geq\tau_0$
\begin{equation}
|\Phi_2|<\delta_0 e^{-\tau}.
\end{equation}
\subsection{A point wise estimate for $\psi$ in $|z|>e^{\frac{\tau}{4}}$}\ \\
The estimate of last section is locally. In this section we make use of comparison principle to obtain an estimate in larger region.
\begin{lemma}
There exists constant $K_2$ such that in the region $|z|>e^{\frac{\tau}{4}}$
\begin{equation}
|\psi|\leq K_2T^{\frac{1}{4}}.
\end{equation}
\end{lemma}
\noindent Proof:\ We separate $\psi$ into 2 parts. Let $\Psi_3$ be the solution satisfies
\begin{equation*}
\begin{cases}
\partial_t\Psi_{3}-\Delta\Psi_3=\lambda_2^{-2}\frac{1}{1+|y_2|^{2+a}}R^{-\frac{1}{2}}\ \ \ \ \ \ \ \ \ \ \ (x,t)\in\mathbb{R}^5\times[0,T).\\
\Psi_3\mid_{t=0}=0,
\end{cases}
\end{equation*}
Denote
\begin{equation*}
\begin{split}
G'=&\lambda_1^{-2}\frac{1}{1+|y_1|^{2+a}}e^{-\tau}R^{-\frac{1}{2}}\textbf{1}_{{|x|}\leq 1}+e^{-\frac{3}{4}\tau}|z|^4\textbf{1}_{e^{\frac{3}{8}\tau}\leq{|z|}\leq 2e^{\frac{3}{8}\tau}}+e^{-\frac{7\tau}{6}}|z|^{4}\textbf{1}_{{|z|}\leq 2e^{\frac{\tau}{2}}}
\\&+e^{-\frac{7\tau}{3}}|z|^{10}\textbf{1}_{{|z|}\leq 2e^{\frac{\tau}{4}}}+e^{-\frac{7\tau}{3}}|z|^{6}\textbf{1}_{{|z|}\leq 2e^{\frac{3}{8}\tau}}+e^{-\frac{7}{3}\tau}\textbf{1}_{{|z|}\leq 1}\\
&+\frac{\delta_0^{\frac{7}{3}}}{1+|x|^{4}}\textbf{1}_{{|z|}\geq e^{\frac{\tau}{4}}}+
\frac{\|Z_2\|_\infty^p}{1+|x|^4}\textbf{1}_{{|z|}\geq e^{\frac{\tau}{2}}}+\frac{e^{-2\tau}}{1+|x|^{3}}\textbf{1}_{{|x|}\geq 1}.
\end{split}
\end{equation*}
Let $\Psi_4$ satisfies
\begin{equation*}
\begin{split}
\begin{cases}
\partial_t\Psi_{4}-\Delta\Psi_4=G'\ \ \ \ \ \ \ \ \ \ \ (x,t)\in\mathbb{R}^5\times[0,T),\\
\Psi_4\mid_{t=0}={(d\cdot e)\eta(\frac{z}{e^{\frac{3}{8}\tau}})}.
\end{cases}
\end{split}
\end{equation*}
Use Lemma 2.4 we obtain
\begin{equation}
|\Psi_3|\leq T^{\frac{1}{4}}.
\end{equation}
Then we check that $\overline{\Psi}=K_2(2e^{-\frac{1}{4}\tau_0}-e^{-\frac{1}{4}\tau})$, $K_2$ sufficient large, is a super solution of $\Psi_4$ in $|z|>e^{\frac{\tau}{4}}$.\\
First, from Lemma 5.3, (5.20) and $$|b_\alpha e_\alpha|\lesssim e^{-\frac{13}{6}\tau}(1+|z|^4),$$ we obtain when $|z|=e^{\frac{\tau}{4}}$
\begin{equation*}|\psi|\lesssim e^{-\tau}.\end{equation*}
Also $\psi=\Psi_3+\Psi_4$, because $|\Psi_3|\leq T^{\frac{1}{4}}$, for $|z|=e^{\frac{\tau}{4}}$
\begin{equation*}|\Psi_4|\leq e^{-\frac{1}{4}\tau_0}\leq\overline{\Psi},\end{equation*}
and for $|z|>e^{\frac{\tau}{4}}$, we verify that\\
\begin{equation*}
\overline{\Psi}_\tau-A_z\overline{\Psi}=\frac{K_2}{4}e^{-\frac{1}{4}\tau}>e^{-\tau}G'>\partial_{\tau}\Psi_4-A_z\Psi_4.
\end{equation*}
For $\tau=\tau_0$
\begin{equation*}
|\Psi_4|\leq e^{-\frac{13}{6}\tau_0}e^{\frac{3}{2}\tau_0}\leq e^{-\frac{2}{3}\tau_0}\leq|\overline{\Psi}|.
\end{equation*}
Use comparison principle $|\psi|\leq\Psi_3+\Psi_4$, we finished the proof.

Now Lemma 5.6 provides an $L^\infty$ estimate, then we give the decay estimate.
\begin{lemma}
For $|x|\geq|q|+1$, there exists constant $K_3>0$ such that
\begin{equation}
|\psi|\leq K_3(\frac{T^3-(T-t)^3}{|x|^3}+\frac{\|Z_2\|_\infty^p}{|x|^2}).
\end{equation}
\end{lemma}
\noindent Proof:\ We are going to check that
$$\widetilde{\psi}=K_3(\frac{T^6-(T-t)^6}{|x|^3}+\frac{\|Z_2\|_\infty^p}{|x|^2}),$$ for $K_3$ sufficient large, is a super solution of $\psi$ in $|x|\geq|q|+1$.
First we verify that if $K_3$ large
\begin{equation*}
\widetilde{\psi}_t-\Delta\widetilde{\psi}=K_3(\frac{3(T-t)^2}{|x|^3}+\frac{2\|Z_2\|_\infty^p}{|x|^4})\geq\psi_t-\Delta\psi.
\end{equation*}
From Lemma 5.6, for $|x|=q+1$
\begin{equation*}
|\psi|\leq K_2T^\frac{1}{4}\leq\widetilde{\psi}.
\end{equation*}
When $t=0$, because the intersection of $|x|\geq|q|+1$ and $$\psi\mid_{t=0}={(d\cdot e)\eta(\frac{z}{e^{\frac{3}{8}\tau}})}$$ is empty.
So in $|x|\geq|q|+1$
\begin{equation*}
\psi \mid_{t=0}=0\leq \frac{\|Z_2\|_\infty^p}{|x|^2}=\widetilde{\psi}\mid_{t=0}.
\end{equation*}
Use comparison principle, we finished the proof.

\section{The proof of theorem 1.1}\label{chap:8}
After the preparation of the last few sections, we start to prove Theorem 1.1.\\
Because for $1\leq i\leq m$
\begin{equation*}
\begin{split}
|b_{\alpha_i} e_{\alpha_i}|\lesssim e^{-\frac{13}{6}\tau}(1+|z|^4),
\end{split}
\end{equation*}
from Lemma 5.3
\begin{equation*}
|\Phi_1(z,\tau)\leq Ce^{-2\tau}(1+|z|^4),\quad \tau\geq\tau_0.
\end{equation*}
combine with Lemma 5.5 and (5.20), we have
\begin{equation}
\begin{split}
|\Phi_2|=&|\int_{\tau_0}^{\tau}e^{A_z(\tau-\tau')}e^{-\tau'}G^\bot d\tau'|\leq \delta_0e^{-\tau}
+e^{-\frac{7}{4}\tau}|z|^4+\frac{e^{-\tau}}{R^{\frac{1}{2}}}|z|^2
\\ &+e^{-6\tau}|z|^{18}+e^{-\frac{10}{3}\tau}|z|^{10}.
\end{split}
\end{equation}
When $|z|\leq  e^{-\frac{\tau}{4}}$
\begin{equation}
\begin{split}
|\Phi_2|&\lesssim \delta_0e^{-\tau}
+e^{-\frac{7}{4}\tau}|z|^4+\frac{e^{-\tau}}{R^{\frac{1}{2}}}|z|^2
+e^{-6\tau}|z|^{18}+e^{-\frac{10}{3}\tau}|z|^{10}\\
&\lesssim \delta_0e^{-\tau}
+e^{-\frac{7}{4}\tau}|z|^4+\frac{e^{-\tau}}{R^{\frac{1}{2}}}|z|^2
+e^{-\frac{5}{2}\tau}|z|^{4}+e^{-\frac{11}{6}\tau}|z|^{4}.\\
\end{split}
\end{equation}
If $\tau_0$ sufficient large, we have
\begin{equation}
|\Phi_2|\leq\delta_0 e^{-\tau}(1+|z|^4),
\end{equation}
and from Lemma 5.6 and Lemma 5.7 we have when $|z|\geq  e^{-\frac{\tau}{4}}$
\begin{equation}
|\psi|\leq\delta_0 \frac{1}{1+|x|^2}.
\end{equation}
We know operator $T(\Phi)=\Phi_*$ maps X in to X, and from local Holder estimate and Ascoli-Arezela theorem we know $T$ is compact. By Schauder fixed point theorem, we know that $T$ has a fixed point $\Phi$ in X. Which finished the proof of Theorem 1.1\\

\bibliographystyle{springer}
\bibliography{mrabbrev,literatur}
\newcommand{\noopsort}[1]{} \newcommand{\printfirst}[2]{#1}
\newcommand{\singleletter}[1]{#1} \newcommand{\switchargs}[2]{#2#1}

\end{document}